\RequirePackage{fix-cm}
\documentclass[smallextended]{svjour3}       
\smartqed  
\usepackage{graphicx}
\usepackage{subfig}
\usepackage{fmtcount}
\usepackage{textcomp}
\usepackage{float}
\usepackage[margin=1in]{geometry}
\usepackage{commath,mathtools}
\usepackage{algorithm}
\usepackage{algorithmic}
\usepackage[english]{babel}
\usepackage[autostyle]{csquotes}
\usepackage{hyperref}
\usepackage{color}
\usepackage{rotating,booktabs}
\usepackage{lscape}
\usepackage[sort,nocompress]{cite}
\usepackage{enumitem}
\usepackage{multirow}
\usepackage{tikz}
\usepackage{comment}
\usetikzlibrary{shapes,arrows,chains}
\usepackage{amssymb}
\usepackage{amsmath}

\usepackage{pifont}

\newcommand{\xbf}{\mathbf{x}}

\newcommand{\zbf}{\mathbf{z}}
\newcommand{\gbf}{\mathbf{g}}
\newcommand{\ubf}{\mathbf{u}}
\newcommand{\vbf}{\mathbf{v}}
\newcommand{\wbf}{\mathbf{w}}

\newcommand{\Wbf}{\mathbf{W}}
\newcommand{\Xbf}{\mathbf{X}}

\newcommand{\etol}{\epsilon_{\mathrm{tol}}}

\def\BibTeX{{\rm B\kern-.05em{\sc i\kern-.025em b}\kern-.08em
    T\kern-.1667em\lower.7ex\hbox{E}\kern-.125emX}}

\begin{document}
\title{A Learned Proximal Alternating Minimization Algorithm and Its Induced Network for a Class of Two-block Nonconvex and Nonsmooth Optimization} 



\author{Yunmei Chen, Lezhi Liu, Lei Zhang
}


\institute{ 
Yunmei Chen,  Lezhi Liu and Lei Zhang\at
              Department of Mathematics, University of Florida, Gainesville, FL 32611, USA \\
              \email{yun, lliu1, leizhang@ufl.edu}           
}

\date{Received: date / Accepted: date}

\maketitle
 
\begin{abstract}
This work proposes a general learned proximal alternating minimization algorithm, LPAM, 
for solving learnable two-block nonsmooth and nonconvex optimization problems.  We tackle the nonsmoothness by an appropriate smoothing technique with automatic diminishing smoothing effect. For smoothed nonconvex problems we modify the proximal alternating  linearized minimization (PALM) scheme by incorporating the residual learning architecture, which has proven to be highly effective in deep network training, and  employing the
block coordinate decent (BCD) iterates as a safeguard for the convergence of the algorithm. We prove that there is a subsequence of the iterates generated by LPAM, which has at least one accumulation point and each accumulation point is a Clarke stationary point. 
 Our method is widely applicable as one can employ various  learning problems formulated as two-block optimizations, and is also easy to be extended for solving multi-block nonsmooth and nonconvex optimization problems. 
The network, whose architecture follows the LPAM exactly,  namely LPAM-net, inherits the convergence properties of the algorithm to make the network  interpretable. As an example application of LPAM-net, we present the numerical and theoretical results on the application of LPAM-net for joint multi-modal MRI reconstruction with significantly under-sampled $k$-space data. The experimental results indicate the proposed LPAM-net is parameter-efficient and has favourable performance in comparison with some state-of-the-art methods.

\end{abstract}




\section{Acknowledgement}
\textbf{Funding} This research was partially supported by NSF grant DMS-2152961 and Simons Foundation Collaboration Grant (Award Number: 584918).\\
\newline
\textbf{Data Availability} The original dataset we used in all experiments are from {\it Multi-modal Brain Tumor Segmentation Challenge 2018}.\cite{menze2014multimodal} 

\section{Introduction}
\label{sec:introduction}

Recent years have witnessed remarkable success of deep learning across various real-world applications. However, a purely data-driven approach may fail to approximate the desired functions, especially when training data are scarce. It is well known that the scarcity of data leads to overfitting and challenges in interpretation. 
To mitigate these issues, the unrolling/unfolding neural networks (UNNs) have been developed and have shown promising results in solving inverse problems arising from computer vision and medical imaging. The UNNs are multi-phase neural networks, in which each phase mimics one iteration of the  optimization algorithms for solving the inverse problems.  
However, despite their promising performance in practice, most of the existing UNNs only superficially resemble steps of optimization algorithms. Consequently, their outputs do not really yield  solutions to  any interpretable variational models. This results is lack of theoretical justifications and convergence guarantees for their outputs.

Recently a novel class of unrolling methods known as learned optimization algorithms (LOA)  \cite{LDA,ldct,wanyuMiccai} has been developed. The goal of LOA is to tackle the challenges associated with solving a class of inverse problems arising from various image reconstruction and synthesis problems. 
These LOAs strategically combine the learned variational model with the prior domain knowledge of underlying physical processes to enhance its interpretability. These methods also utilize optimization to guarantee the convergence.

The optimization scheme induces a highly structured deep network, aligning its architecture precisely with the iterative algorithm.  Hence, the network inherits  the convergence property of the algorithm, and it outputs an approximation of the solution to the variational model. As a result, it is interpretable and parameter efficient. However, the existing LOAs only consider a single block of variables.

Motivated by a wide range of  applications in machine learning that involve multi-modalities or multi-domains problems,  such as  multi-task learning, transfer learning, multi-modal learning and fusion, dual domain image reconstruction, and    image synthesis,  in this work we extend the idea of LDA to deal with  learned  multi-block nonconvex and nonsmooth optimization problems. More specifically, we develop a novel convergent learned proximal alternating minimization (LPAM) algorithm to solve the following  class of learnable  optimization problems:
 \begin{equation} \label{M1}
      (M_1):\quad \mbox{minimize} \quad \Phi(\xbf_{1}, \xbf_{2}; \Theta):=H_1(\xbf_{1}; \theta_1)+H_2(\xbf_{2}; \theta_2)+H(\xbf_{1}, \xbf_{2}; \theta), \quad  (\xbf_{1}, \xbf_{2})\in \mathbb R^n\times \mathbb R^m,
\end{equation}
 where each learnable function is a parameterized function. $\Theta=(\theta_1, \theta_2, \theta)$ is the set of  learned parameters. Each function in \eqref{M1} is possibly nonconvex and nonsmooth. We shall also study the convergence and iteration complexity of the algorithm, and provide experimental results on the application on the joint multi-contrast MRI  reconstruction with significantly under-sampled data.


\section{Related work}
\label{sec:related work}

\subsection{ Related work on alternating minimization algorithms:}

Alternating minimization (AM) or “two block” coordinate minimization (BCM) algorithm plays an indispensable role in solving a  rich class of problems formulated as follows:
\begin{equation} \label{M}
    (M):\quad \mbox{minimize} \quad \Phi(\xbf_{1}, \xbf_{2}):=H_1(\xbf_{1})+H_2(\xbf_{2})+H(\xbf_{1}, \xbf_{2}), \quad  (\xbf_{1}, \xbf_{2})\in \mathbb R^n\times \mathbb R^m. 
\end{equation}
The standard approach to solve this problem is the so-called Gauss-Seidel iteration scheme, also known as block coordinate descent (BCD) method. This scheme alternately keeps one block of updated variable fixed and optimize the other block \cite{palomar2010convex, bertsekas2015parallel}. 
With assumptions on convexity and smoothness, several convergence results have been established \cite{beck2013convergence,charbonnier1997deterministic}.

Inspired by many practical problems arising  from machine learning and image processing, such as non-negative matrix factorization,  blind image deconvolution, low rank matrix completion \cite{jain2013low,hardt2014understanding}, phase retrieval problem \cite{netrapalli2013phase}, interference alignment\cite{peters2009interference}, and compressed sensing using synthesis sparsity models\cite{abolghasemi2012gradient}
joint multi-modal image reconstruction\cite{sroubek2011robust}, the AM algorithms  for nonconvex (including biconvex) and nonsmooth minimization have attracted great interest. One of the effective approaches to solve nonconvex and nonsmooth  problems \eqref{M} is the proximal alternating minimization (PAM) algorithms. The PAM algorithm developed in \cite{attouch2010proximal} can be viewed as a proximal regularization of a two block Gauss-Seidel method for solving \eqref{M}. Under the assumption that  the objective function satisfies Kurdyka-       Lojasiewicz (KL) property, the convergence to a critical point of $\Phi$ of each bounded sequence generated by the algorithm has been obtained in \cite{attouch2013convergence}. 
However, PAM algorithm requires exact minimization of a nonconvex and
nonsmooth problem in each step. To overcome this difficulty, the pioneer work  \cite{bolte2014proximal} introduced  a proximal alternating linearized minimization (PALM) algorithm,  which  can be viewed as alternating the steps of the proximal forward-backward scheme. Building on the KL property, it was proved  that each bounded sequence generated by PALM globally converges to a critical point. Later, the work \cite{pock2016inertial} proposed an inertial version of the PALM (iPALM) algorithm motivated from the Heavy Ball method of Polyak \cite{polyak1964some}. It is also closely related to multi-step algorithms of Nesterov's accelerated gradient method \cite{nesterov2018lectures}.
The global convergence to a critical point is also obtained for iPALM under the assumption that the objective function $\Phi$ satisfies the KL property and the same smoothness conditions as that for PALM.

\subsection{Related work on unrolling neural networks inspired by AM algorithms:}

In recent years, the unrolling neural networks inspired by AM algorithms unroll the iteration scheme of an AM algorithm  into a multi-phase deep neural network. 
Several unrolling convolutional neural networks (CNN) based on AM algorithms have been developed for solving inverse problems in computer vision including
the unrolling neural networks resembling  the steps of alternating direction method of multipliers (ADMM) algorithm, and block coodinate decent (BCD) algorithm. 
The ADMM-net \cite{sun2016deep}  with the architecture derived from  ADMM has been applied to compressive sensing  MRI reconstruction. 
The deeply aggregated alternating minimization
 neural network (DeepAM)  \cite{kim2017deeply} for image restoration uses a data-driven approach in the energy minimization framework. At each iteration, the DeepAM
advances  the following two steps in the conventional AM algorithm. One step uses a data-driven approach to update the auxiliary variable representing the gradient of the image by a CNN that plays the role of regularization. The other step recovers the image by minimizing the data fidelity and the disparity between the gradient of the image and  updated auxiliary variable. 
Moreover, several  versions of BCD-nets that resemble the steps of BCD algorithm have been successfully applied to  solve the inverse problems for signal/image recovery/reconstruction   \cite{aggarwal2018modl, chun2019bcd, chun2018deep, chun2020convolutional,chun2020momentum}.

The common idea of those BCD-nets is  combining a denoising CNN in defining the regularization as one module, namely a denoising module, and a model based image reconstruction (MBIR) as the other module into the BCD framework. The MBIR module enhances the data fidelity for the image from  the denoising module. 
For instance, at each iteration of the BCD-net proposed in \cite{chun2018deep}, it trains an image mapping CNN  using identical convolutional kernels in both encoders and decoders as one module in a BCD framework for image recovery. The other module in this BCD-net minimizes an objective function consisting of the data fidelity and  penalty between to be updated image and the deniosed image from image mapping CNN. 
This BCD-Net \cite{chun2018deep} is modified from \cite{chun2019bcd}, where a sparsity promoting denoising network is trained in denoising modules and the accelerated proximal gradient method using a majorizer is applied to MBIR modules with statistical CT data-fidelity  for low-dose CT reconstruction.
Furthermore, to achieve fast and convergent solutions for the inverse problem in image processing, the work in \cite{chun2020momentum} proposes a Momentum-Net. The Momentum-Net is motivated by applying the block proximal extrapolated gradient method using a Majorizer and convolutional autoencoders to MBIR.  Each iteration of Momentum-Net consists of three core modules: image refining, extrapolation and MBIR. The image refining module is a denoising module that trains an image refining CNN to remove iteration-wise artifact. The extrapolation module uses momentum from previous updates to amplify the changes in subsequent iterations and accelerate convergence. The MBIR is non-iterative, since the MBIR problem in consideration of this work   has a closed-form solution. Momentum-Net guarantees convergence to a fixed-point for general differentiable
nonconvex data-fit terms and convex feasible sets, under the conditions that the sequence of paired refining CNNs is asymptotically nonexpansive and the refined image from the image refining module gives an asymptotically block-coordinate minimizer.

Very recently, \cite{ding2023learned} presented a learned alternating algorithm aimed at solving a specific variational model for dual-domain sparse-view CT reconstruction. In this model, the joint term is convex and smooth, and the regularization terms are learnable in both the image and sinogram domains, respectively. In this study, we extend their work by considering a broader scenario, where each term in the objective function can be nonconvex and nonsmooth. We provide sufficient conditions for taking smoothing approximations of nonsmooth objective functions. We also provide the analysis of the iteration complexity in addition to the convergence analysis.


\section{Contributions and Organization}
In this work we propose a novel learned proximal alternating minimization algorithm addressing nonsmooth and nonconvex two-block variational models with provable convergence and iteration complexity. 
The proposed LPAM algorithm for solving ($M_1$)
determines the architecture of the deep neural network, hence the design of the
LPAM considers  not only  the convergence and efficiency, but also the ability to assist the training of the network parameters $\Theta$, i.e., reducing the error in minimizing the loss function.

Our main idea for designing this algorithm is as follows: (i) We tackle the nonsmoothness by an appropriate smoothing technique with automatic diminishing smoothing effect; (ii) We modify the PALM scheme by incorporating the residual learning architecture, which has proven to be highly effective in deep network training \cite{he2016deep}; and (iii) When some modified PALM iterates  fail to satisfy certain conditions, we employ the block coordinate decent (BCD) iterates as a safeguard to ensure convergence. Moreover, we prove that a subsequence generated by the proposed LPAM has accumulation points and all of them are Clarke stationary points of the nonsmooth nonconvex problem.This algorithm can be readily extended to multi-block nonsmooth and nonconvex optimization.

This algorithm significantly broadens its applicability, as it relaxes the strict convergence conditions of the existing algorithms.
The PALM or iPALM algorithm assumes the joint term $H(\xbf_{1}, \xbf_{2})$ to be smooth and possibly nonconvex, and the proximal points associated with $H_i(\xbf_{i})$ ($i=1,2$) are easy to obtain. Their global convergence results are  built on KL property of the objective function. These assumptions in general cannot be met by the learnable parameterized functions  representing neural networks. 
 The proposed LPAM algorithm will not require these conditions.
 Each function in \eqref{M1} could be nonconvex and nonsmooth. 
  Without those restrictions, the LPAM algorithm is  flexible and applicable to deep learning problems in various applications. Moreover, the convergence of LPAM can still be established in the sense of  sub-sequence convergence to Clarke stationary points (see detail in Section 3). But it cannot have global convergence to a critical point as the PALM or iPALM algorithm.

The convergence result of the BCD-net is based on the assumptions that the paired image refining neural networks are asymptotically nonexpansive,
 the data fidelity is differentiable, and the MBIR problem is not difficult to solve. 
The proposed LPAM algorithm does not require these assumptions. 
The convergence of LPAM is assured by the convergence of the BCD algorithm for smooth optimization and the property of the limit of the gradient of the smoothed objective functions as smoothing factor tends zero. By removing the constraints on the model, the PALM can significantly broaden its range of applications. 

The remainder of the paper is organized as follows. In Section 5  we construct the LPAM algorithm motivated by the PALM algorithm and the BCD algorithm, and introduce the three steps that constitute the LPAM algorithm. In Section 6, we prove the convergence of the LPAM algorithm and deduce the iteration complexity of this algorithm. In Section 7, we apply the LPAM algorithm to the joint reconstruction of T1 and T2 MRI images with significantly under-sampled data to show the promising performance of the proposed method. Furthermore, we provide convergence analysis and comparison with several existing methods. All the figures are listed at the end of the article.

\bigskip
Throughout the paper, we use the following notations without further notification. 
\begin{enumerate}
    \item $H(\xbf_{1},\xbf_{2};\theta)$ denotes the function $H$ of $\xbf_{1}$ and $\xbf_{2}$ given a collection of parameters $\theta$.  $H_{\epsilon}(\xbf_{1},\xbf_{2};\theta)$ denotes the smooth approximation of $H(\xbf_{1},\xbf_{2};\theta)$ with a smoothing parameter $\epsilon$.
    \item  For a differentiable function $f(\xbf_1,\xbf_2)$, we denote the gradient of $f$ with respect to $\xbf_{i}$  by $\nabla_i f(\xbf_1,\xbf_2)$  $(i=1,2)$. The gradient $f$ on the domain $(\xbf_1,\xbf_2)\in R^n \times R^m$ is represented by $\nabla_{1,2}f(\xbf_1,\xbf_2)$.
    \item  We use $\|\cdot\|$ to denote the standard $L^2$ norm in Euclidean space, 
    and write $\langle \cdot , \cdot \rangle _{\mathbb{R}^{n}}$ for the standard inner product on $\mathbb{R}^{n}$.
    \item $\chi:=\mathbb R^n\times \mathbb R^m$. 
\end{enumerate} 

\section{Learned Proximal Alternating Minimization (LPAM) Algorithm }

In this section we propose a learned proximal alternating minimization (LPAM) algorithm  for solving the problem \eqref{M1}.
Since $\theta$ is learned, we will omit it and simply write \eqref{M1}  as \eqref{M}. 

For $H_1$,$H_2$, $H$ and $\Phi$ we assume
\begin{align*}
&(A1): H_1, \ H_2 \ \mbox{and} \ H \ \mbox{ are proper and  locally Lipschitz on} \ \mathbb R^n, \ \mathbb R^m \ \mbox{and} \ \mathbb R^n \times \mathbb R^m, \mbox{respectively}. \\
&\quad \mbox{Each of them is possibly non-convex and non-smooth.}\\
&(A2): \Phi \mbox{ is coercive}, 
\mbox{and} \  \min_{\Xbf \in \chi} \Phi(\Xbf) > -\infty
\end{align*}



The proposed LPAM algorithm consists of three stages. In the first stage we convert the nonconvex and nonsmooth problem \eqref{M} to a nonconvex smooth optimization problem by an appropriate smoothing procedure. The smoothing method is not unique, but we require the smooth approximation  $\Phi_{\epsilon}$ of $\Phi:=H_1(\xbf_1)+H_2(\xbf_2)+H(\xbf_1, \xbf_2)$, which is defined as
\begin{equation} \label{smoothed phi}
    \Phi_{\epsilon}(\xbf_1, \xbf_2):=H_{1,\epsilon}(\xbf_1)+H_{2,\epsilon}(\xbf_2)+H_{\epsilon}(\xbf_1, \xbf_2)
\end{equation}
to satisfy following conditions:

\begin{itemize}
 \item ({\bf C1}:) $\nabla H_{1,\epsilon}$, $\nabla H_{2,\epsilon}$ and $\nabla_{1,2}H_{\epsilon}$ are Lipschitz continuous in $\mathbb R^n$, $\mathbb R^m$ and $\mathbb R^n\times \mathbb R^m$, respectively. 
 
    \item ({\bf C2}:)  $\lim_{\epsilon\to 0^{+},\zbf_i\to \xbf_i}H_{i,\epsilon}(\zbf_i)=H_i(\xbf_i)$, 
    $i=1,2$, for any given $\xbf_i$ in the corresponding domain. Also for any $\Xbf=(\xbf_1,\xbf_2)\in \mathbb R^n\times \mathbb R^m$, we assume that when $\epsilon\to 0^{+}$ and ${\bf Z}=(\zbf_1,\zbf_2)\to \Xbf$,
    $\lim_{\epsilon\to 0^{+},{\bf Z}\to \Xbf}H_{\epsilon}({\bf Z})=H(\Xbf)$.
     \item ({\bf C3}:) There is a continuous and non-negative function $\mathfrak{m}: [0,\infty)\to [0,\infty)$ such that $\mathfrak{m}(0)=0$,
     $$\Phi_{\epsilon}(\xbf_1,\xbf_2)+\mathfrak{m}(\epsilon)\le \Phi_{\delta}(\xbf_1,\xbf_2)+\mathfrak{m}(\delta)\,\, \mbox {for all} \,(\xbf_1,\xbf_2)\in \chi, \, \mbox {and}\, 
     0<\epsilon\le \delta.$$ 
    \item ({\bf C4}:) For any $\Xbf^*=(\xbf_1, \xbf_2)\in\mathbb R^n\times \mathbb R^m$ and any $\Xbf^k\to \Xbf^*$,  $\lim_{\Xbf^k\to \Xbf^*,\epsilon\to 0^{+}}\nabla_{1,2}\Phi_{\epsilon}(\Xbf^k)\in \partial^c\Phi(\Xbf^*)$, where $\partial^c\Phi(\Xbf^*)$ stands for the Clarke subdifferential of $\Phi$ (see Definition \ref{def-clark-gen}).
\end{itemize}

\begin{remark} The Lipschitz constants of $\nabla H_{1,\epsilon},\nabla H_{2,\epsilon}$ and $\nabla_{1,2}H_{\epsilon}$ may tend to infinity as $\epsilon$ tends to $0$. The condition ${\bf C3}$ is much weaker than the monotonicity of $\Phi_{\epsilon}$. The condition ${\bf C3}$ roughly says that as long as $\Phi_{\epsilon}$ is not too much different from a monotonic function with respect to $\epsilon$, our method still leads to convergence. This brings great convenience in application. 
\end{remark}

In the second stage we solve the smoothed nonconvex problem with a fixed smoothing factor $\epsilon$, i.e. 
\begin{equation}
 \label{SM}   
\min_{\xbf_1, \xbf_2}\{\Phi_{\epsilon}(\xbf_1, \xbf_2):=H_{1,\epsilon}(\xbf_1)+H_{2,\epsilon}(\xbf_2)+H_{\epsilon}(\xbf_1, \xbf_2)\}.
\end{equation}
In light of the substantial improvement in practical performance by ResNet \cite{he2016deep},  we propose a modified  PALM that incorporates the architecture of the ResNet for solving \eqref{SM}.  With $\epsilon = \epsilon_k >0$, PALM \cite{bolte2014proximal} generates a sequence of iterates $\{\xbf_1^{k+1},\xbf_2^{k+1}\}$ by

$$\zbf_1^{k+1}=\xbf_1^k-\alpha_k \nabla H_{1,\epsilon_k}(\xbf_1^k),$$
$$
  \ubf_1^{k+1}=\arg\min_{\ubf} \frac{1}{2\alpha_k}\|\ubf-\zbf_1^{k+1}\|^2+H_{\epsilon_k}(\ubf,\xbf_2^k), 
$$
and
$$\zbf_2^{k+1}=\xbf_2^k-\beta_k\nabla H_{2,\epsilon_k}(\xbf_2^k),$$
$$
  \ubf_2^{k+1}=\arg\min_{\ubf} \frac{1}{2\beta_k}\|\ubf-\zbf_2^{k+1}\|^2+H_{\epsilon_k}(\ubf_1^{k+1}, \ubf),  
$$
where $\alpha_k$ and $\beta_k$ are step sizes. In some situations, the proximal points $\ubf_1^{k+1}$ and $\ubf_2^{k+1}$ are not easy to compute, so we replace $H_{\epsilon_k}(\ubf,\xbf_2^k)$  by its linear approximation at $\zbf_1^{k+1}$ together with a proximal term $\frac{1}{2p_k} \| \ubf-\zbf_1^{k+1} \|^2$. Similarly, we replace $H_{\epsilon_k}(\ubf_1^{k+1}, \ubf)$  by its linear approximation at $\zbf_2^{k+1}$ together with a proximal term $\frac{1}{2q_k} \| \ubf-\zbf_2^{k+1} \|^2$. Thus, we have
\begin{align} \label{u1update}
    \ubf_1^{k+1} &=\arg\min_{\ubf} \frac{1}{2\alpha_k}\|\ubf-\zbf_1^{k+1}\|^2+<\nabla_1H_{\epsilon_k}(\zbf_1^{k+1},\xbf_2^k), \ubf-\zbf_1^{k+1}>+\frac{1}{2p_k} \| \ubf-\zbf_1^{k+1} \|^2 \\ \nonumber
   & =\zbf_1^{k+1}-\tau_k\nabla_1H_{\epsilon_k}(\zbf_1^{k+1},\xbf_2^k), 
 \quad  \quad \mbox{with} \quad \tau_k = \frac{\alpha_k p_k}{\alpha_k + p_k}.\nonumber
 \end{align}
 where we recall that $\nabla_1H_{\epsilon_k}$ stands for the gradient for the first $n$ variables.
 \begin{align}\label{u2update}
 \ubf_2^{k+1} &=\arg\min_{\ubf} \frac{1}{2\beta_k}\|\ubf-\zbf_2^{k+1}\|^2+<\nabla_2H_{\epsilon_k}(\ubf_1^{k+1}, \zbf_2^{k+1}), \ubf-\zbf_2^{k+1}>+\frac{1}{2q_k} \| \ubf-\zbf_2^{k+1} \|^2 \\ 
 &=\zbf_2^{k+1}-\gamma_k\nabla_2H_{\epsilon_k}(\ubf_1^{k+1}, \zbf_2^{k+1}),  \quad  \quad \mbox{with} \quad \gamma_k = \frac{\beta_k q_k}{\beta_k + q_k}. \nonumber
 \end{align}
 In deep learning approach, the step sizes $\alpha_k$, $\tau_k$, $\beta_k$ and $\gamma_k$ could be learnable hyper-parameters. The advantage for taking the scheme \eqref{u1update} and \eqref{u2update} to update $u_1$ and $u_2$ is that this scheme is consistent with the architecture of residual learning \cite{he-1,he-2},  which only learns the correction needed for $\zbf_i^{k+1}$ $(i=1,2)$ and avoids gradient vanishing in network training, and thus improves solution quality in practice.  However, the convergence of the sequence $\{(\ubf_1^{k+1}, \ubf_2^{k+1})\}$ is not guaranteed. Inspired by the proof of convergence in \cite{bolte2014proximal}, we propose the following: 

If $\ubf_1^{k+1}$ and $\ubf_2^{k+1}$ satisfy
\begin{equation}\label{gen-condtion1}
    \Phi_{\epsilon_k}(\ubf_1^{k+1},\ubf_2^{k+1})-\Phi_{\epsilon_k}(\xbf_1^k,\xbf_2^k)\le -a(\|\ubf_1^{k+1}-\xbf_1^k\|^2+\|\ubf_2^{k+1}-\xbf_2^k\|^2) 
\end{equation}
and 
\begin{equation}\label{gen-condtion2}
\|\nabla_{1,2} \Phi_{\epsilon_k}(\xbf_1^k,\xbf_2^{k})\|\le a^{-1}(\|\ubf_1^{k+1}-\xbf_1^k\|+\|\ubf_2^{k+1}-\xbf_2^k\|), 
\end{equation}
for some $a>0$, we take $\xbf_1^{k+1}=\ubf_1^{k+1}$, $\xbf_2^{k+1}=\ubf_2^{k+1}$.

If one of (\ref{gen-condtion1}) and (\ref{gen-condtion2}) is violated, we compute $(\vbf_1^{k+1},\vbf_2^{k+1})$ by BCD algorithm with a simple line-search strategy for better step sizes as follows: Let $\bar \alpha,\bar \beta\in (0,1)$, we compute
\begin{align} \label{v1v2}
   & \vbf_1^{k+1}=\arg\min_{\vbf}<\nabla_1 \Phi_{\epsilon_k}(\xbf_1^k,\xbf_2^k), \vbf-\xbf_1^k>+\frac 1{2\bar \alpha}\|\vbf-\xbf_1^k\|^2, \\
& \vbf_2^{k+1}=\arg\min_{\vbf}<\nabla_2 \Phi_{\epsilon_k}(\vbf_1^{k+1},\xbf_2^k), \vbf-\xbf_2^k>+\frac 1{2\bar \beta}\|\vbf-\xbf_2^k\|^2. 
\end{align}
The first order optimality condition leads to
\begin{align}\label{v-12}
&\vbf_1^{k+1}=\xbf_1^k-\bar \alpha( \nabla H_{1,\epsilon_k}(\xbf_1^k)+\nabla_1H_{\epsilon_k}(\xbf_1^k,\xbf_2^k)), 
\\
&\vbf_2^{k+1}=\xbf_2^k-\bar \beta( \nabla H_{2,\epsilon_k}(\xbf_2^k)+\nabla_2 H_{\epsilon_k}(\vbf_1^{k+1},\xbf_2^k)).\nonumber
\end{align}
If for some $\delta\in (0,1)$,
\begin{equation}\label{v-condition-6}
\Phi_{\epsilon}(\vbf_1^{k+1},\vbf_2^{k+1})-\Phi_{\epsilon}(\xbf_1^k,\xbf_2^k)
\le -\delta (\|\vbf_1^{k+1}-\xbf_1^k\|^2+\|\vbf_2^{k+1}-\xbf_2^k\|^2), 
\end{equation}
we set
$$\xbf_1^{k+1}=\vbf_1^{k+1},\quad \xbf_2^{k+1}=\vbf_2^{k+1}, $$
otherwise we reduce step sizes $(\bar \alpha,\bar \beta)$ by
 $\rho (\bar \alpha,\bar \beta)\rightarrow (\bar \alpha,\bar \beta)$ where $0<\rho<1$, and recompute  $\vbf_1^{k+1},\vbf_2^{k+1}$ until the condition (\ref{v-condition-6}) holds. Then take $\xbf_1^{k+1}=\vbf_1^{k+1}$, $\xbf_2^{k+1}=\vbf_2^{k+1}$. 

In the third stage, we check if the $\ell_{2}$-norm of the gradient of the smoothed objective function has been reduced enough, so that we can  perform the second stage with a reduced smoothing factor. By gradually reducing the smoothing factor, we obtain a subsequence of the iterates that converges to a Clarke stationary point of the original nonconvex and nonsmooth problem.

The scheme of LPMA is given below (For notational simplicity, we omit all learned network parameters).

\begin{algorithm}[H]
\caption{
Learned Proximal Alternating Minimization (LPAM) Algorithm }
\label{alg:2}
\begin{algorithmic}[1]
\STATE \textbf{Input:} Initial $\xbf_1^0$, $\xbf_2^0$, $0<\rho, \gamma, \delta<1$, and $ \bar \alpha, \bar \beta, \sigma, a, \epsilon_0>0$. Maximum iteration $K$ or tolerance $\varepsilon_{tol}>0$. 
\FOR{$k=0,1,2,\dots,K$} 
\STATE $\zbf_1^{k+1}=\xbf_1^k-\alpha_k\nabla H_{1,\epsilon_k}(\xbf_1^k)$,
\STATE $\ubf_1^{k+1}=\zbf_1^{k+1}-\tau_k \nabla_1 H_{\epsilon_k}( \zbf_1^{k+1}, \xbf_2^k)$, 
\STATE $\zbf_2^{k+1}=\xbf_2^k-\beta_k\nabla  H_{2,\epsilon_k}(\xbf_2^k)$
\STATE $\ubf_2^{k+1}=\zbf_2^{k+1}-\gamma_k \nabla_2 H_{\epsilon_k}(\ubf_1^{k+1}, \zbf_2^{k+1}, )$ 
\IF{ conditions \eqref{gen-condtion1} and \eqref{gen-condtion2} hold,} 
\STATE set $\xbf_1^{k+1}=\ubf_1^{k+1}$, $\xbf_2^{k+1}=\ubf_2^{k+1}$,
\ELSE  
\STATE For given $\bar \alpha,\bar \beta\in (0,1)$
\STATE $\vbf_1^{k+1}=\xbf_1^{k}-\bar \alpha (\nabla H_{1,\epsilon_k}(\xbf_1^k)+\nabla_1 H_{\epsilon_k}(\xbf_1^k,\xbf_2^k)), $ \label{marker}
\STATE $\vbf_2^{k+1}=\xbf_2^{k}-\bar \beta(\nabla  H_{2,\epsilon_k}(\xbf_2^k)+\nabla_2 H_{\epsilon_k}(\vbf_1^{k+1},\xbf_2^k)).$
\IF{ condition \eqref{v-condition-6} holds,}
\STATE set $\xbf_1^{k+1}=\vbf_1^{k+1}$, $\xbf_2^{k+1}=\vbf_2^{k+1}$,
\ELSE
\STATE update $(\bar \alpha, \bar \beta) \leftarrow \rho (\bar \alpha, \bar \beta)$,
then \textbf{go to}~\ref{marker},
\ENDIF
\ENDIF
\STATE \textbf{if} $\|\nabla_{1,2}\Phi_{\varepsilon_k}(\xbf_1^{k+1}, \xbf_2^{k+1})\| < \sigma \gamma {\varepsilon_k}$,  set $\varepsilon_{k+1}= \gamma {\varepsilon_k}$;  \textbf{otherwise}, set $\varepsilon_{k+1}={\varepsilon_k}$.
\STATE \textbf{if} $\sigma {\varepsilon_k} < \varepsilon_{tol}$, terminate.
\ENDFOR
\STATE \textbf{Output:} $\xbf_1^{k+1}$, $\xbf_2^{k+1}$
\end{algorithmic}
\end{algorithm}

We would like to comment here for the LPAM algorithm. Line 2 to line 18 can be viewed as the steps of the inner iteration, which, for a fixed $\varepsilon$, generates a sequence $\{\xbf_1^k,\xbf_2^k\}$ such that $\|\nabla_{1,2}\Phi_{\epsilon}(\xbf_1^k,\xbf_2^k)\|\to 0$ as $k\to\infty$ (see the first statement of Lemma 1 in the next section). This  guarantees the reduction of $\|\nabla_{1,2}\Phi_{\varepsilon}(\xbf_1^k,\xbf_2^k)\|$ for a fixed $\varepsilon$. When this value drops below a prefixed value, line 19  decreases $\varepsilon$, and then we repeat the procedures from line 2 to line 18 to generates a new sequence $\{\xbf_1^k,\xbf_2^k\}$ corresponding to the reduced $\varepsilon$. We continue these steps  until the termination condition in the line 20 is met. The sequence formed by the iterates, in which the reduction criterion for $\varepsilon$ is met,  has at least one accumulation point and each accumulation point is a Clarke Stationary Point of the original problem. This will be proved in Theorem 1 in Section 6.

{\bf Remark:} In order to make the LPMA algorithm match the ResNet architecture, if the function $H_{\epsilon}(\xbf_1, \xbf_2)$ is learnable and $H_{1,\epsilon}(\xbf_1)$ and $H_{2,\epsilon}(\xbf_2)$ are given, it is better to use the scheme above. If the functions $H_{1,\epsilon}(\xbf_1)$ and $H_{2,\epsilon}(\xbf_2)$ are learnable and $H_{\epsilon}(\xbf_1, \xbf_2)$ is given, it is better to change the steps 3-6 to the following:
\begin{align}
& \zbf_1^{k+1}=\xbf_1^k-\alpha_k \nabla_1 H_{\epsilon_k}(\xbf_1^{k}, \xbf_2^k), \\
& \ubf_1^{k+1}=\zbf_1^{k+1}-\tau_k \nabla H_{1,\epsilon_k}(\zbf_1^{k+1}); \\
& \zbf_2^{k+1}=\xbf_2^k-\beta_k \nabla_2 H_{\epsilon_k}(\ubf_1^{k+1}, \xbf_2^{k} ), \\
& \ubf_2^{k+1}=\zbf_2^{k+1}-\gamma_k \nabla H_{2,\epsilon_k}(\zbf_2^{k+1});
\end{align}
and keep the other steps in Algorithm \ref{alg:2} unchanged. 

\section{Convergence Analysis} \label{sec: method}
Since we deal with a nonconvex and nonsmooth optimization problem,  the Clarke subdifferential, which is  based on the concept of generalized directional derivatives, is employed here 
to characterize the optimality of the solutions to \eqref{M}.

\begin{definition} \label{def-clark-gen} (Clarke subdifferential). Suppose that $f:\mathbb R^n\times \mathbb R^m\to (-\infty,\infty]$ is locally Lipschitz. The Clarke subdifferential of $f$ at $\Xbf=(\xbf_1,\xbf_2)$ 
is defined as 
\begin{align*}
\partial^c f(\Xbf):=\{\Wbf\in \mathbb R^n\times \mathbb R^m | <\Wbf, \textbf{V}>\le \limsup_{\textbf{Z} \to \Xbf,t\downarrow 0}\frac{f(\textbf{Z}+t\textbf{V})-f(\textbf{Z})}t,\quad \forall \; \textbf{V} \in \mathbb R^n\times \mathbb R^m\}.
\end{align*}

where $\Xbf,\textbf{Z}\in \mathbb R^n\times \mathbb R^m$, $<\Wbf,\textbf{V}>=<\wbf_1,\vbf_1>_{\mathbb R^n}+<\wbf_2,\vbf_2>_{\mathbb R^m}$ for $\textbf{V}=(\vbf_1,\vbf_2)$, $\Wbf=(\wbf_1,\wbf_2)$.
\end{definition}

\begin{definition}\label{def-1} (Clarke stationary point) For a locally Lipschitz function $f$ defined as in Definition \ref{def-clark-gen}, a point $\Xbf=(\xbf_1, \xbf_2)\in \mathbb R^n\times \mathbb R^m$ is called a Clarke Stationary Point of $f$, if $0\in \partial^c f(\textbf{X})$. 
\end{definition}


\begin{lemma}\label{intro-lem-1}
Assume that $\Phi_{\epsilon}(\xbf_1, \xbf_2)$ is a smooth approximation of $\Phi(\xbf_1, \xbf_2)$ defined in \eqref{smoothed phi} satisfying the conditions ({\bf C1})-({\bf C4}). Let $\varepsilon, \eta, a >0$, $\rho, \bar \alpha,\bar \beta \in (0,1)$, and $\Xbf^0=(\xbf_1^0,\xbf_2^0)$ be an arbitrary initial condition. Suppose $\{\Xbf^k=(\xbf_1^k,\xbf_2^k)\}$ is the sequence generated by repeating Lines 2--18 in Algorithm 1 with $\varepsilon_k = \varepsilon$. Then  we have 
\begin{enumerate}
    \item $\|\nabla_{1,2}\Phi_{\epsilon}(\xbf_1^k,\xbf_2^k)\|\to 0$ as $k\to\infty$.
    
    \item  Let $(\bar \alpha,\bar \beta)$ $(\bar \alpha,\bar \beta\in (0,1))$ be the initial step size for $(\vbf_1^k,\vbf_2^k)$. The maximum number of required line search steps $\ell_{max}$ is 
 $$\ell_{max}=\Bigg[\frac{\log \big((\frac{L_{\epsilon}}2+\delta)\max\{\bar \alpha,\bar \beta\}\big )}{\log 1/\rho}\Bigg]+1, $$
    where $\big[A\big]$ represents the largest integer less than or equal to $A$, and $L_{\epsilon}$ is the sum of the Lipschitz constants for $\nabla_1 H_{1,\epsilon}$, $\nabla_2 H_{2,\epsilon}$, and $\nabla_{1,2} H_{\epsilon}$, i.e.
    $$L_{\epsilon}:=L_{\nabla_1 H_{1,\epsilon}}+L_{\nabla_2 H_{2,\epsilon}}+L_{\nabla_{1,2} H_{\epsilon}}.$$ 
    \item For any $\eta>0$,
\begin{equation} \label{mink}
    \min\{k\in \mathbb N;\quad \|\nabla_{1,2} \Phi_{\epsilon}(\xbf_1^{k+1},\xbf_2^{k+1})\|<\eta \}\le (2 a^{-3}+\frac{4\max(\bar \alpha,\bar \beta)^2}{\delta \min(\bar \alpha,\bar \beta)^2\rho^2}L_{\epsilon}^2)\frac{\Phi(\Xbf^0)-\Phi_*+1}{\eta^2}\end{equation}
\end{enumerate} 
 
\end{lemma}

\noindent{\bf Proof of Lemma \ref{intro-lem-1}:} 

Given $(\xbf_1^k, \xbf_2^k)$, in the case that $(\ubf_1^{k+1},\ubf_2^{k+1})$ generated by Algorithm 1 satisfies the conditions \eqref{gen-condtion1}-\eqref{gen-condtion2}, we take $(\xbf_1^{k+1},\xbf_2^{k+1})=(\ubf_1^{k+1},\ubf_2^{k+1})$. From \eqref{gen-condtion1}-\eqref{gen-condtion2}, it is easy to get
\begin{equation}\label{grad-b-1}\|\nabla_{1,2} \Phi_{\epsilon}(\xbf_1^{k},\xbf_2^{k})\|^2\le \frac{2(\Phi_{\epsilon}(\xbf_1^{k},\xbf_2^{k})-\Phi_{\epsilon}(\ubf_1^{k+1},\ubf_2^{k+1}))}{a^3},
\end{equation}


If $(\ubf_1^{k+1},\ubf_2^{k+1})$ fails to satisfy  the condition \eqref{gen-condtion1} or \eqref{gen-condtion2}, the algorithm  computes $(\vbf_1^{k+1},\vbf_2^{k+1})$ by  \eqref{v-12}  with a line search strategy. Let $\ell_k$
be the number of the required line search steps to meet the condition \eqref{v-condition-6}, and then  we take  $(\xbf_1^{k+1},\xbf_2^{k+1})=(\vbf_1^{k+1},\vbf_2^{k+1})$. This gives the following
\begin{equation}\label{v-12c}
\vbf_1^{k+1}=\xbf_1^k-\bar \alpha \rho^{\ell_k}\nabla_1\Phi_{\epsilon}(\xbf_1^k,\xbf_2^k), \qquad \vbf_2^{k+1}=\xbf_2^k-\bar \beta \rho^{\ell_k}\nabla_2\Phi_{\epsilon}(\vbf_1^{k+1},\xbf_2^k).
\end{equation}

Since $\nabla_{1,2} \Phi_{\epsilon}$ is $L_{\epsilon}$-Lipschitz continuous, we have
\begin{align}
\Phi_{\epsilon}(\vbf_1^{k+1},\vbf_2^{k+1})&
\le \Phi_{\epsilon}(\vbf_1^{k+1},\xbf_2^k)+\nabla_2\Phi_{\epsilon}(\vbf_1^{k+1},\xbf_2^k)\cdot
(\vbf_2^{k+1}-\xbf_2^k)+\frac{L_{\epsilon}}{2}\|\vbf_2^{k+1}-\xbf_2^k\|^2\nonumber\\
\le &\Phi_{\epsilon}(\xbf_1^k,\xbf_2^k)+\nabla_1\Phi_{\epsilon}(\xbf_1^k,\xbf_2^k)\cdot (\vbf_1^{k+1}-\xbf_1^k)+\frac{L_{\epsilon}}{2}\|\vbf_1^{k+1}-\xbf_1^k\|^2\nonumber\\
+&\nabla_2\Phi_{\epsilon}(\vbf_1^{k+1},\xbf_2^k)\cdot
(\vbf_2^{k+1}-\xbf_2^k)+\frac{L_{\epsilon}}{2}\|\vbf_2^{k+1}-\xbf_2^k\|^2 \label{compare-phi-3}.
\end{align}
The combination of \eqref{v-12c} and \eqref{compare-phi-3} yields
\begin{equation}\label{v-decay}
\Phi_{\epsilon}(\vbf_1^{k+1},\vbf_2^{k+1})\leq   \Phi_{\epsilon}(\xbf_1^k,\xbf_2^k)+(-\frac 1{\bar \alpha\rho^{\ell_k}}+\frac{L_{\epsilon}}{2})\|\vbf_1^{k+1}-\xbf_1^k\|^2
+(-\frac{1}{\bar \beta\rho^{\ell_k}}+\frac{L_{\epsilon}}{2})\|\vbf_2^{k+1}-\xbf_2^k\|^2. 
\end{equation}
 Hence the condition (\ref{v-condition-6}) is met if 
 \begin{equation}\label{add-1}
-\frac 1{\bar \alpha\rho^{\ell_k}}+\frac{L_{\epsilon}}{2}\le -\delta, \quad \mbox{and} \quad
-\frac{1}{\bar \beta\rho^{\ell_k}}+\frac{L_{\epsilon}}{2}\le -\delta. 
 \end{equation}
Hence, for any $k=1,2,\ldots$ the maximum line search steps $\ell_{max}$ required for \eqref{v-condition-6} satisfies 
$$\rho^{\ell_{max}}= (\delta+L_{\varepsilon}/2)^{-1}(\max\{\bar \alpha,\bar \beta\})^{-1},$$ so we abuse the notation  $\ell_{max}$ by defining it as
 \begin{equation}\label{line-bound}
 \ell_{max} =\Bigg[\frac{\log \big((\frac{L_{\epsilon}}2+\delta)\max\{\bar \alpha,\bar \beta\}\big )}{\log 1/\rho}\Bigg]+1
\end{equation}
This proves the second statement of this lemma. Moreover, the fact that $\ell_k\leq \ell_{max}$ for any $k=0,1,2,\ldots $, implies
     \begin{equation}\label{stepsize bound}
    \rho^{\ell_{k}}\geq  \rho^{\ell_{max}}=
    \frac{\rho}{ (\delta+L_{\varepsilon}/2)\max\{\bar \alpha,\bar \beta\}}.
    \end{equation}

Next we prove the other two statements. Note that from (\ref{v-12c})  the condition (\ref{v-condition-6}) can be rewritten as: 
 \begin{align}\label{line-up-2}
\Phi_{\epsilon}(\vbf_1^{k+1},\vbf_2^{k+1})\le \Phi_{\epsilon}(\xbf_1^k,\xbf_2^k)
-\delta\bar \alpha^2\rho^{2\ell_k}\|\nabla_1\Phi_{\epsilon}(\xbf_1^k,\xbf_2^k)\|^2
-\delta \bar \beta^2 \rho^{2\ell_k}
\|\nabla_2\Phi_{\epsilon}(\vbf_1^{k+1},\xbf_2^k)\|^2 .
\end{align}
Now we estimate the last term in \eqref{line-up-2} in terms of $\nabla_{1,2}\Phi_{\epsilon}(\xbf_1^{k},\xbf_2^k)$.
First by the Lipschitz continuity of $\nabla_2\Phi_{\epsilon}$  we have 
$$\|\nabla_2 \Phi_{\epsilon}(\vbf_1^{k+1},\xbf_2^k)-\nabla_2\Phi_{\epsilon}(\xbf_1^k,\xbf_2^k)\|\le L_{\epsilon}\|\vbf_1^{k+1}-\xbf_1^k\|.$$

Clearly this implies
\begin{equation}\label{ineq-c}\|\nabla_2 \Phi_{\epsilon}(\vbf_1^{k+1},\xbf_2^k)\|\ge \|\nabla_2 \Phi_{\epsilon}(\xbf_1^k,\xbf_2^k)\|-L_{\epsilon}\|\vbf_1^{k+1}-\xbf_1^k\|.
\end{equation}
Here if the right hand side is positive, we will choose  $\sigma\in (0,1)$ to be determined, and take advantage
of the following inequality:
$(a-b)^2\ge \frac{\sigma}2a^2-\sigma b^2$ together with \eqref{v-12c}
to obtain
\begin{align}\label{tem-1}
  & \|\nabla_2\Phi_{\epsilon}(\vbf_1^{k+1},\xbf_2^k)\|^2 \nonumber\\
  \ge &\frac{\sigma}2 \|\nabla_2\Phi_{\epsilon}(\xbf_1^k,\xbf_2^k)\|^2-\sigma L_{\epsilon}^2\|\vbf_1^{k+1}-\xbf_1^k\|^2   \nonumber \\
 = & \frac{\sigma}2\|\nabla_2\Phi_{\epsilon}(\xbf_1^k,\xbf_2^k)\|^2-\sigma L_{\epsilon}^2\bar \alpha^2\rho^{2\ell_k}\|\nabla_1\Phi_{\epsilon}(\xbf_1^k,\xbf_2^k)\|^2.
\end{align}
Inserting (\ref{tem-1}) into \eqref{line-up-2}, we have
\begin{align} \label{line-up-3}
&\Phi_{\epsilon}(\vbf_1^{k+1},\vbf_2^{k+1})-\Phi_{\epsilon}(\xbf_1^k,\xbf_2^k)\\
  \le &-\big (\delta \bar \alpha^2 \rho^{2\ell_k}-\delta \bar \beta^2 \rho^{2\ell_k}L_{\epsilon}^2\sigma\bar \alpha^2\rho^{2\ell_k}\big )\|\nabla_1\Phi_{\epsilon}(\xbf_1^{k},\xbf_2^k)\|^2 
  -\delta \bar \beta^2\rho^{2\ell_k}\frac{\sigma}2\|\nabla_2\Phi_{\epsilon}(\xbf_1^{k},\xbf_2^k)\|^2 \nonumber \\
  =&-\delta \bar \alpha^2\rho^{2\ell_k}(1-\bar \beta^2\rho^{2\ell_k}L_{\epsilon}^2\sigma)\|\nabla_1\Phi_{\epsilon}(\xbf_1^{k},\xbf_2^k)\|^2 -\delta \bar \beta^2\rho^{2\ell_k}\frac{\sigma}2\|\nabla_2\Phi_{\epsilon}(\xbf_1^{k},\xbf_2^k)\|^2 \nonumber 
\end{align}
If $\bar \beta^2\rho^{2\ell_k}L_{\epsilon}^2<\frac 23$, we just choose $\sigma=\frac 12$. Otherwise we set $\sigma$ as
$\sigma=1/(2\bar \beta^2\rho^{2\ell_k}L_{\epsilon}^2)$.
In either case the coefficient of $\|\nabla_1 \phi_{\epsilon}(\xbf_1^k,\xbf_2^k)\|^2$ is less than $-\frac{\delta}2\bar \alpha^2\rho^{2\ell_k}$. For the coefficient of $\|\nabla_2\phi_{\epsilon}(\xbf_1^k,\xbf_2^k)\|^2$, in the first case it is less than $-\frac{\delta}4\bar \beta^2\rho^{2\ell_k}$. In the second case it is less than $-\frac{\delta}{4L_{\epsilon}^2}$. Then it is easy to verify that in both cases, by (\ref{stepsize bound})  both coefficients are less than $-\frac{\delta \min(\bar \alpha,\bar \beta)^2 \rho^2}{4\max (\bar \alpha,\bar \beta)^2 L_{\epsilon}^2}$
Thus we have 
\begin{equation}\label{line-up-n}\Phi_{\epsilon}(\vbf_1^{k+1},\vbf_2^{k+1})-\Phi_{\epsilon}(\xbf_1^k,\xbf_2^k)
\le -D_{\epsilon} \|\nabla_{1,2}\Phi_{\epsilon}(\xbf_1^k,\xbf_2^k)\|^2,
\end{equation}
where 
\begin{equation}\label{d-ep}
D_{\epsilon}=\frac{\delta \min(\bar \alpha,\bar \beta)^2\rho^2}{4\max(\bar \alpha,\bar \beta)^2L_{\epsilon}^2}.
\end{equation}
Consequently,
\begin{equation}\label{grad in v}
  \|\nabla_{1,2}\Phi_{\epsilon}(\xbf_1^{k},\xbf_2^k)\|^2 \leq \frac{4\max(\bar \alpha,\bar \beta)^2L_{\epsilon}^2}{\delta \min(\bar \alpha,\bar \beta)^2\rho^2}\{
\Phi_{\epsilon}(\xbf_1^k,\xbf_2^k)-\Phi_{\epsilon}(\vbf_1^{k+1},\vbf_2^{k+1})\}.
\end{equation}
If the right hand side of (\ref{ineq-c}) is negative, we obtain from (\ref{v-12c}) that 
\[\|\nabla_2\Phi_{\epsilon}(x_1^k,x_2^k)\|\le L_{\epsilon}\|v_1^{k+1}-x_1^k\|=L_{\epsilon}\bar \alpha\rho^{l_k}\|\nabla_1\Phi_{\epsilon}(x_1^k,x_2^k)\|.\]
Consequently,
\begin{equation}\label{line-up-6}
\|\nabla_{1,2}\Phi_{\epsilon}(x_1^k,x_2^k)\|^2=\|\nabla_1 \Phi_{\epsilon}(x_1^k,x_2^k)\|^2+\|\nabla_2 \Phi_{\epsilon}(x_1^k,x_2^k)\|^2\le (1+L_{\epsilon}^2\bar \alpha^2\rho^{2 l_k})\|\nabla_1 \Phi_{\epsilon}(x_1^k,x_2^k)\|^2.
\end{equation}
Then we use (\ref{line-up-6}) to estimate (\ref{line-up-2}) differently:
\begin{align}\label{line-up-5}
&\Phi_{\epsilon}(\vbf_1^{k+1},\vbf_2^{k+1})- \Phi_{\epsilon}(\xbf_1^k,\xbf_2^k)\\
\le &-\delta\bar \alpha^2\rho^{2\ell_k}\|\nabla_1\Phi_{\epsilon}(\xbf_1^k,\xbf_2^k)\|^2
\le -\frac{\delta \bar \alpha^2\rho^{2l_k}}{1+L_{\epsilon}^2\bar \alpha^2\rho^{2l_k}}\|\nabla_{1,2}\Phi_{\epsilon}(x_1^k,x_2^k)\|^2.\nonumber
\end{align}
If $L_{\epsilon}^2\bar \alpha^2\rho^{2l_k}\ge 1$, we have
\[\Phi_{\epsilon}(\vbf_1^{k+1},\vbf_2^{k+1})- \Phi_{\epsilon}(\xbf_1^k,\xbf_2^k)\le -\frac{\delta}{2L_{\epsilon}^2}\|\nabla_{1,2}\Phi_{\epsilon}(x_1^k,x_2^k)\|^2\]
Otherwise 
\[\Phi_{\epsilon}(\vbf_1^{k+1},\vbf_2^{k+1})- \Phi_{\epsilon}(\xbf_1^k,\xbf_2^k)\le -\frac{\delta}2\bar \alpha^2\rho^{2l_k}\|\nabla_{1,2}\Phi_{\epsilon}(x_1^k,x_2^k)\|^2\]
In either case we see that (\ref{line-up-n}) still holds. 
Note that we used (\ref{stepsize bound}) for a lower bound of $\rho^{2l_k}$. Thus we have established (\ref{grad in v}).

Observing \eqref{gen-condtion1}, \eqref{grad-b-1}, \eqref{v-condition-6} and \eqref{grad in v} we see that in either case of $(\xbf_1^{k+1},\xbf_2^{k+1})=(\ubf_1^{k+1},\ubf_2^{k+1})$ or $(\xbf_1^{k+1},\xbf_2^{k+1})=(\vbf_1^{k+1},\vbf_2^{k+1})$, there are two constants 
\begin{equation}\label{two-cons}
b_1=\min\{a, \delta \}\quad \mbox{and} \quad b_2=\max\{2a^{-3},   \frac{4\max(\bar \alpha,\bar \beta)^2L_{\epsilon}^2}{\delta \min(\bar \alpha,\bar \beta)^2\rho^2}  \}, 
\end{equation}
such that 
\begin{equation}\label{x-decay}
    \Phi_{\epsilon}(\xbf_1^{k+1},\xbf_2^{k+1})-\Phi_{\epsilon}(\xbf_1^k,\xbf_2^k)
\le -b_1 (\|\xbf_1^{k+1}-\xbf_1^k\|^2+\|\xbf_2^{k+1}-\xbf_2^k\|^2), 
\end{equation}
and 
\begin{equation}\label{x-grad}
   \|\nabla_{1,2}\Phi_{\epsilon}(\xbf_1^{k},\xbf_2^k)\|^2 \leq b_2(
\Phi_{\epsilon}(\xbf_1^k,\xbf_2^k)-\Phi_{\epsilon}(\xbf_1^{k+1},\xbf_2^{k+1})). 
\end{equation}

From \eqref{x-decay} and {\bf C2} we get 
\begin{equation} \label{phi-eps-decay}
 \Phi_{\epsilon}(\xbf_1^{k+1},\xbf_2^{k+1})<\Phi_{\epsilon}(\xbf_1^k,\xbf_2^k)< \ldots < \Phi_{\epsilon}(\xbf_1^0,\xbf_2^0)
\leq \Phi(\xbf_1^0,\xbf_2^0)+1
\end{equation}


Furthermore, for any integer $K>0$ summing up \eqref{x-grad} for $k = 0, \ldots, K$, we have 
\begin{equation}\label{add-e-2}
\sum_{k=0}^K\|\nabla_{1,2}\Phi_{\epsilon}(\xbf_1^{k},\xbf_2^k)\|^2\le b_2(\Phi_{\epsilon}(\xbf_1^{0},\xbf_2^0)-\Phi_{\epsilon}(\xbf_1^{K+1},\xbf_2^{K+1}))\le b_2(\Phi_{\epsilon}(\xbf_1^{0},\xbf_2^0)-\Phi_*+1).
\end{equation}
We use $\Phi_*$ to denote a uniform lower bound for all $\Phi_{\epsilon}$: $\min_{\xbf}\Phi_{\epsilon}(\xbf)\ge \Phi_*$ for all $\epsilon>0$. Since $b_2$ is independent of $K$, we see that for any fixed $\epsilon>0$,
$\|\nabla_{1,2}\Phi_{\epsilon}(\xbf_1^{k},\xbf_2^k)\|^2\to 0$ as $k\to \infty$.
This proves the first statement of Lemma \ref{intro-lem-1}.
For the third statement of Lemma \ref{intro-lem-1}, we set $\kappa:=\min\{k\in\mathbb N;\,\, \|\nabla_{1,2}\Phi_{\epsilon}(\xbf_1^{k+1},\xbf_2^{k+1})\|\le \eta\}$, then we know that $\|\Phi_{\epsilon}(\xbf_1^{k+1},\xbf_2^{k+1})\|\ge \eta$ for all $k\le \kappa-1$. Thus from (\ref{add-e-2}) we have
$$\kappa \eta^2\le \sum_{k=0}^{\kappa-1}\|\nabla_{1,2} \Phi_{\epsilon}(\xbf_1^{k+1},\xbf_2^{k+1})\|^2=\sum_{k=1}^{\kappa}\|\nabla_{1,2} \Phi_{\epsilon}(\xbf_1^k,\xbf_2^k)\|^2\le b_2(\Phi_{\epsilon}(\xbf_1^0,\xbf_2^0)-\Phi_*+1).$$
The third statement follows immediately.
Lemma \ref{intro-lem-1} is established. $\Box$

\medskip

From the first statement of Lemma 1 the reduction criterion for $\epsilon_k$ in Line 19 of Algorithm 1 must be satisfied within finitely many iterations for any $k=1,2,\ldots$. Hence $\epsilon_k$ will eventually be small enough to terminate the algorithm. 
Let $k_{l} $ be the counter of iteration when the criterion in Line 19 of Algorithm 1 is met for the $l$-th time (we set $k_{0}=-1$), then we can partition the iteration counters $k=0,1,2,\dots,$ into segments accordingly, so that in the $l$-th segment $k=k_{l} +1,\dots,k_{l+1} $  and $\epsilon_k = \epsilon_{k_{l} +1} = \epsilon_0 \gamma^l$. In  the next theorem we will give the length of each segment, which will lead to the iteration complexity  of Algorithm \ref{alg:2} for any $\etol>0$.

\begin{theorem}\label{convergence-thm} Let $\Phi_{\epsilon}$ and $\Phi$ be described as in Lemma \ref{intro-lem-1} and ${\bf (C1)-(C4)}$ hold for $\Phi_{\epsilon}$. If $\epsilon_{tol}=0$,
suppose $\{\Xbf^k=(\xbf_1^k,\xbf_2^k)\}$ is the sequence generated by Algorithm 1 with arbitrary initial condition $\Xbf^0=(\xbf_1^0,\xbf_2^0)$. Let
$\{\tilde \Xbf^l\}=\{(\tilde \xbf_1^l,\tilde \xbf_2^l)=:(\xbf_1^{k_l+1},\xbf_2^{k_l+1})\}$ be the subsequence of $\{\Xbf^k\}$ where the reduction criterion for $\epsilon$ in the algorithm is met for $k=k_l$ and $l=1,2...$. Then $\{\tilde \Xbf^l\}$ has at least one accumulation point and each accumulation point is a Clarke Stationary Point of the original problem. Moreover, the number of iterations, $k_{l+1}-k_l$, for the $l-th$ segment is bounded by 
\begin{equation}\label{inner-c}
k_{l+1}-k_l\le (2a^{-3}+\frac{4\min\{\bar \alpha,\bar \beta\}^2}{\delta \max\{\bar \alpha,\bar \beta\}^2}L_{\epsilon_0\gamma^l}^2)\frac{\Phi(\Xbf^0)-\Phi_*+1}{\sigma^2\epsilon_0^2\gamma^{2l+2}}.
\end{equation}
For any $\epsilon_{tol}>0$, the total number of iterations for Algorithm 1 to terminate with $\epsilon_{tol}$ is bounded by
\begin{equation}\label{add-e-5}
\sum_{l=1}^{l_0-1}(2a^{-3}+\frac{4\min\{\bar \alpha,\bar \beta\}^2}{\delta \max\{\bar \alpha,\bar \beta\}^2}L_{\epsilon_0\gamma^l}^2)\frac{\Phi(\Xbf^0)-\Phi_*+1}{\sigma^2\epsilon_0^2\gamma^{2l+2}}=O(L_{\epsilon_{tol}}^2\epsilon_{tol}^{-2})
\end{equation}
where $l_0$ is number of reductions and we have  $l_0-1=\frac{\log \sigma\epsilon_0/\epsilon_{tol}}{\log 1/\gamma}$, $C_1>0$ depends only on $\rho$ and $\delta$.
\end{theorem}

\begin{remark} In application we usually have $L_{\epsilon}\sim \epsilon^{-1}$, which makes the bound in (\ref{add-e-5}) $O(\epsilon_{tol}^{-4})$.
\end{remark}

\begin{remark} The figure below shows the generation of the sequence $\{\tilde \Xbf^l\}$ in Theorem 1. The iterates highlighted in red indicate the selected  sequence $\{\tilde \Xbf^l\}$, where the reduction criterion for $\varepsilon_k$ is met at $k=k_l$ for $l=0,1,2...$.
\begin{figure}[H]
\centering
\includegraphics[width=0.70\textwidth]{xep.pdf}
\label{xep}
\end{figure}
\end{remark}
\noindent{\bf Proof of Theorem \ref{convergence-thm}:} 
 First we claim that the sequence $\tilde \Xbf^l=(\tilde \xbf_1^l,\tilde \xbf_2^l)$ is compact. This is largely based on {\bf C3}. In the first step, the initial $\Xbf^0$ corresponds to $\epsilon_0$. Then we have $\Phi_{\epsilon_0}(\Xbf^1)< \Phi_{\epsilon_0}(\Xbf^0)$. After this we may have $\Xbf^2,..,\Xbf^l$ and then we have $\epsilon_1=\gamma \epsilon_0$ and 
 $\Phi_{\epsilon_1}(\Xbf^{l+1})<\Phi_{\epsilon_1}(\Xbf^l)$. Here we mention that we don't have 
 $\Phi_{\epsilon_1}(\Xbf^{l+1})\le \Phi_{\epsilon_0}(\Xbf^l)$, but by the third requirement of $\Phi_{\epsilon}$ we have 
 $$\Phi_{\epsilon_1}(\Xbf^{l+1})+\mathfrak{m}(\epsilon_1)\le \Phi_{\epsilon_0}(\Xbf^{l+1})+\mathfrak{m}(\epsilon_0)
 <\Phi_{\epsilon_0}(\Xbf^0)+\mathfrak{m}(\epsilon_0). $$
 In a similar manner we then prove that $\Phi_{\epsilon_l}(\tilde \Xbf^l)+\mathfrak{m}(\epsilon_l)$ is uniformly bounded for all $l$. Since $\mathfrak{m}$ is a positive continuous function and all $\epsilon_l\le \epsilon_0$, after removing the bound for all $\mathfrak{m}(\epsilon_l)$ we have the uniform bound for all 
 $\Phi_{\epsilon_l}(\tilde \Xbf^l)$.

Next we estimate the iteration times based on the previous lemma. If $\epsilon_{tol}\neq 0$, to estimate $k_{l+1}-k_l$, which is the iteration times for the inner circle, we see that $\epsilon=\epsilon_{k_l}=\epsilon_0 \gamma^l$, then $\eta$ in Lemma \ref{intro-lem-1} is $\eta=\sigma \gamma \epsilon_{k_l+1}=\sigma\epsilon_0\gamma^{l+1}$, the initial step is $\Xbf_{k_l+1}$. Then
$$k_{l+1}-k_l\le (2a^{-3}+\frac{4\min\{\bar \alpha,\bar \beta\}^2}{\delta \max\{\bar \alpha,\bar \beta\}^2}L_{\epsilon_0\gamma^l}^2)\frac{\Phi(\Xbf^0)-\Phi_*+1}{\sigma^2\epsilon_0^2\gamma^{2l+2}}. $$
and (\ref{inner-c}) is justified.
Let $l_0$ be the number of reductions. Then $\sigma \epsilon_0\gamma^{l_0-1}\ge \epsilon_{tol}$. Thus 
\begin{equation}\label{est-l-0}
l_0-1\le \frac{\log \sigma \epsilon_0/\epsilon_{tol}}{\log 1/\gamma } .
\end{equation}
Based on this we have
$$\sum_{l=0}^{l_0-1}(k_{l+1}-k_l)\le \sum_{l=1}^{l_0-1}(2a^{-3}+\frac{4\min\{\bar \alpha,\bar \beta\}^2}{\delta \max\{\bar \alpha,\bar \beta\}^2}L_{\epsilon_0\gamma^l}^2)\frac{\Phi(\Xbf^0)-\Phi_*+1}{\sigma^2\epsilon_0^2\gamma^{2l+2}},$$
which yields the left hand side of  (\ref{add-e-5}). If we use $L_{\epsilon_{tol}}$ as an upper bound of $L_{\epsilon_0\gamma^l}$ for all $l$, the summation above leads to $O(L_{\epsilon_{tol}}^2\gamma^{2 l_0})$. Using the estimate of $l_0$ in (\ref{est-l-0}) we see that this bound is $O(L_{\epsilon_{tol}}^2\epsilon_{tol}^{-2})$. 
Since $\|\nabla_{1,2}\Phi_{\epsilon_{k_l}}(\Xbf^{k_l+1})\|\le \sigma \gamma \epsilon_{k_l}=\sigma \epsilon_0 \gamma^{l+1} \to 0$ as $l\to \infty$, so if we choose a subsequence of 
$\{\xbf_1^{k_l+1},\xbf_2^{k_l+1}\}$ that converges to an equilibrium point $\bar \xbf=(\bar \xbf_1,\bar \xbf_2)$, by {\bf C4}, $0\in \partial^c\phi$, then Definition \ref{def-1} says $\bar \xbf$ is a Clarke Stationary Point of $\Phi$. 
Theorem \ref{convergence-thm} is established. $\Box$

\medskip

\subsection{Small initial steps} Finally in this subsection we prove that if the initial steps $\bar \alpha$ and $\bar \beta$ satisfy
\begin{equation}\label{step-final}
\delta<\bar \alpha L_{\epsilon}<1,\quad \delta<\bar \beta L_{\epsilon}<1, 
\end{equation}
for some $\delta>0$, we can improve the estimate of $D_{\epsilon}$ from (\ref{d-ep}) to $D_{\epsilon}\le cL_{\epsilon}^{-1}$ for some universal constant $c>0$. 

Since under this assumption we have $-1+\frac{L_{\epsilon}}2\bar \alpha <0$ and $-1+\frac{L_{\epsilon}}2\bar \beta<0$, we can estimate the decrease of $\Phi_{\epsilon}(\xbf_1^{k+1},\xbf_2^{k+1})$ in terms of $\|\nabla_{1,2} \Phi_{\epsilon}(\xbf_1^k,\xbf_2^k)\|^2$
using \eqref{v-12c},\eqref{v-decay} and $\ell_k=0$:

\begin{align}\label{a-new-1}
\Phi_{\epsilon}(\vbf_1^{k+1},\vbf_2^{k+1})\le \Phi_{\epsilon}(\xbf_1^k,\xbf_2^k)
+(-\frac{1}{\bar \alpha }+\frac{L_{\epsilon}}2)(\bar \alpha)^{2} \|\nabla_1 \Phi_{\epsilon}(\xbf_1^k,\xbf_2^k)\|^2\nonumber\\
+(-\frac{1}{\bar \beta }+\frac{L_{\epsilon}}2)(\bar \beta)^2\|\nabla_2\Phi_{\epsilon}(\vbf_1^{k+1},\xbf_2^k)\|^2.
\end{align}

To estimate $\|\nabla_{1,2} \Phi_{\epsilon}(\vbf_1^{k+1},\xbf_2^k)\|$ we use Lipschitz continuity and an elementary inequality to obtain (\ref{tem-1}) as before. Inserting (\ref{tem-1}) into (\ref{a-new-1}) we have

\begin{align}\label{a-new-3}
&\Phi_{\epsilon}(\vbf_1^{k+1},\vbf_2^{k+1})-\Phi_{\epsilon}(\xbf_1^k,\xbf_2^k)\nonumber\\
\le &\bigg ( (-\frac{1}{\bar \alpha }+\frac{L_{\epsilon}}2)(\bar \alpha )^2-(-\frac{1}{\bar \beta}+\frac{L_{\epsilon}}2)(\bar \beta )^2\sigma L_{\epsilon}^2(\bar \alpha )^2\bigg )\|\nabla_1 \Phi_{\epsilon}(\xbf_1^k,\xbf_2^k)\|^2\nonumber\\
&\quad +\frac{\sigma}2(-\frac{1}{\bar \beta}+\frac{L_{\epsilon}}2)(\bar \beta )^2\|\nabla_2\Phi_{\epsilon}(\xbf_1^k,\xbf_2^k)\|^2\nonumber\\ 
=&c_{1,\epsilon}\|\nabla_1\Phi_{\epsilon}(\xbf_1^k,\xbf_2^k)\|^2+c_{2,\epsilon}\|\nabla_2\Phi_{\epsilon}(\xbf_1^k,\xbf_2^k)\|^2.
\end{align}
where
\begin{align}\label{c-1e}
c_{1,\epsilon}=(-\frac{1}{\bar \alpha }+\frac{L_{\epsilon}}2)(\bar \alpha )^2-(-\frac{1}{\bar \beta }+\frac{L_{\epsilon}}2)(\bar \beta )^2\sigma 
L_{\epsilon}^2(\bar \alpha )^2\nonumber\\
=\bar \alpha \bigg ((-1+\frac{L_{\epsilon}}2\bar \alpha )-(-1+\frac{L_{\epsilon}}2\bar \beta )\sigma (L_{\epsilon}\bar \beta )(L_{\epsilon}\bar \alpha )\bigg )
\end{align}
$$c_{2,\epsilon}=\frac{\sigma}2(-\frac{1}{\bar \beta }+\frac{L_{\epsilon}}2)(\bar \beta )^2
=\frac{\sigma}2(-1+\frac{\bar \beta L_{\epsilon}}2)(\bar \beta ).$$

we choose $\sigma$ to be 
$$\sigma=\frac{1-\frac{L_{\epsilon}}2\bar \alpha }{2(L_{\epsilon}\bar \beta )(L_{\epsilon}\bar \alpha )(1-\frac{L_{\epsilon}}2\bar \beta )}. $$

Then 
$$c_{1,\epsilon}=\frac{\bar \alpha}2(-1+\frac{L_{\epsilon}\bar \alpha}2),
\quad c_{2,\epsilon}=\frac{-1+L_{\epsilon}/2 \bar \alpha}{4L_{\epsilon}(L_{\epsilon}\bar \alpha)}. $$
Observing (\ref{step-final}) we have 
$$|c_{1,\epsilon}|+|c_{2,\epsilon}|\le CL_{\epsilon} $$
for some $c>0$ independent of $\epsilon$. This clearly implies that $D_{\epsilon}$ is bounded by $CL_{\epsilon}^{-1}$ for some universal $C>0$.

\section{Experiments} 

In this section, we examine the performance of the proposed method for joint MRI T1 and T2 image reconstruction using significantly under-sampled data. We outline the proposed variational model, explore the convergence property of the LPAM algorithm used to solve the model, and demonstrate the efficiency of the LPAM-net. 


The content is organized into five subsections. The first subsection introduces the variational model, especially the learned joint feature extractor. Then in the second subsection we verify the assumptions for the convergence of the algorithm. Then in subsection 3, we describe the dataset used in the experiments and the metrics used to evaluate the reconstruction results. Subsection 4 elucidates the complete reconstruction procedure including the initialization network and the LPAM-net. Finally, in the last subsection, we demonstrate the numerical results of the experiments via comparisons with several existing methods.


\subsection{Proposed model and the feature extractor}
\label{experiment}
Given the under-sampled k-space data ${\textbf{f}_{1},\textbf{f}_{2}}$ of the 
modalities T1 and T2, our goal is to jointly reconstruct the corresponding images $(\xbf_{1}$ and $\xbf_{2})$. We formulate the T1-T2 joint reconstruction as the following nonconvex and nonsmooth minimization problem:
\begin{equation} \label{modell}
\min_{\Xbf=(\xbf_1,\xbf_2)^T}\frac 12 \|\textbf{P} F\xbf_1-\textbf{f}_1\|^2+\frac 12\|\textbf{P}F\xbf_2-\textbf{f}_2\|^2+\|g_{\theta}(\xbf_1,\xbf_2)\|_{2,1}. 
  \end{equation}
where $\xbf_1\in\mathbb R^n$, $\xbf_2\in\mathbb R^n$, $\textbf{f}_1\in\mathbb C^n$ and $\textbf{f}_2\in\mathbb C^n$. $\textbf{f}_1$ and $\textbf{f}_2$ are the  k-space data, and $\gbf$ is a learned joint feature extractor mapping from $\mathbb C^{n\times 2}$ to $\mathbb C^{n\times d},$ where $n$ is the spatial dimension and $d$ is the channel number of the output feature tensor. Here, $F$ stands for the discrete Fourier transform and $\textbf{P}$ is the binary matrix representing the k-space mask when acquiring data for $\xbf_{1}$ and $\xbf_{2}.$ $(\xbf_1,\xbf_2)$ means the image $\xbf_1$ and the image $\xbf_2$ are concatenated. In this model, the first two terms are data fidelity terms and the last term is the regularization term. The regularization term $\gbf_{\theta}$ is a learnable common feature extractor for $\xbf_{1}$ and $\xbf_{2}$ and it enhances the sparsity of the common features. 

More specifically, let $\Xbf=(\xbf_1,\xbf_2)^T$  denote the vector of the concatenated two images of dimension $2n$. 
Then, the regularization term is represented by $\|\gbf_{\theta}(\Xbf)\|_{2,1}\ $.
The $\gbf_{\theta,i}(\Xbf)\in \mathbb C^{d}$ is the $i$-feature vector of $\Xbf$ for $i=1,..., n$. Here, $n$ represents the pixel number of the image. $\nabla_{1,2} \gbf_{\theta}(\Xbf)$ represents the gradient of $\gbf_{\theta}$ as a function of $2n$ variables.  Then,
$$\|\gbf_{\theta}(\Xbf)\|_{2,1}=\sum_{i=1}^n\|\gbf_{\theta, i}(\Xbf)\|, $$
where $\gbf$ is a vanilla l-layer CNN with component-wise activation function $\sigma$. 
\begin{equation}\label{g}
\gbf_{\theta}(\Xbf)=\wbf_l(\sigma...\sigma(\wbf_3*\sigma(\wbf_2*\sigma*(\wbf_1* (\Xbf) )))
\end{equation}
Here $\sigma$ is a smoothed ReLu function defined as:

\begin{equation}\label{sigma}
\sigma(x)= \left\{\begin{array}{ll}
0,\quad &\mbox{if}\quad x\le -\delta\\
\frac{1}{4\delta}x^2+\frac 12 x+\frac{\delta}4,\quad &\mbox{if}\quad -\delta<x<\delta\\
x,\quad &\mbox{if}\quad x\ge \delta. 
\end{array}
\right. 
\end{equation}

We denote $\gbf_{\theta, i}(\Xbf)$ as $\gbf_{i}(\Xbf)$ for simplicity in the following text and in the convergence analysis.

Let $r_{\epsilon}(\Xbf)$ be the smoothed approximation of the regularization term with the parameter $\epsilon$ as follows:
$$r_{\epsilon}(\Xbf)=\sum_{i\in I_0}\frac{1}{2\epsilon}\|\gbf_i(\Xbf)\|^2+\sum_{i\in I_1}(\|\gbf_i(\Xbf)\|-\frac{\epsilon}2).  $$
where 
$I_0=\{i\in [n] | \|\gbf_i(\Xbf)\|\le \epsilon \}$, $I_1=[n]\setminus I_0$. Then, the gradient of the smoothed term is:

\begin{equation}\label{deltar}
\nabla_{1,2} r_{\epsilon}(\Xbf)=\sum_{i\in I_0}\nabla_{1,2} \gbf_i(\Xbf)^{\top}\frac{\gbf_i(\Xbf)}{\epsilon}+\sum_{i\in I_1}\nabla_{1,2} \gbf_i(\Xbf)^{\top}\frac{\gbf_i(\Xbf)}{\|\gbf_i(\Xbf)\|}. 
\end{equation}
Now the smoothed objective function $\Phi_{\epsilon}(\Xbf)$ of $\Phi(\Xbf)$ is written as
$$\Phi_{\epsilon}(\Xbf)=\frac 12 \|\textbf{P} F\xbf_1-\textbf{f}_1\|^2+\frac 12\|\textbf{P}F\xbf_2-\textbf{f}_2\|^2+r_{\epsilon}(\Xbf).$$

\subsection{Convergence Analysis}
\subsubsection{Theoretical proof}
In this subsection, we verify that the proposed variational model \eqref{modell} and its smoothing approximation 
satisfy the assumptions required in Lemma 1 and Theorem 1. This verification provides a strong theoretical support for the convergence of T1-T2 joint reconstruction. 

{\bf Step one: Verification of ${\bf C1}$  and ${\bf C2}$}: 

From the definition of $r_{\epsilon}$, it is easy to see that both ${\bf C1}$ and ${\bf C2}$ hold.

{\bf Step two: Verification of ${\bf C3}$}:

In order to prove {\bf C3}, We choose the continuous function $\mathfrak{m}(\epsilon)=\dfrac{n}{2}\epsilon$ and $\mathfrak{m}(\delta)=\dfrac{n}{2}\delta$. Then, we need to prove the following inequality:
\begin{equation}\label{mono-1}\Phi_{\epsilon}(\Xbf)+\dfrac{n}{2}\epsilon\le \Phi_{\delta}(\Xbf)+\dfrac{n}{2}\delta\,\, \mbox {for all} \, \Xbf \in \chi, \,
 \mbox{and}\, \, 0<\epsilon\le \delta.
 \end{equation}

Note that the data fidelity term is independent of smoothing parameter, hence we only need to prove the following inequality:

\begin{equation}\label{third-a}
r_{\epsilon,i}(\Xbf)+\frac{\epsilon}{2}\le r_{\delta,i}(\Xbf)+\frac{\delta}{2} \quad, 
 \mbox{for all} \quad  \Xbf \in \chi, \quad \mbox{and} \quad 0<\epsilon<\delta.
\end{equation}


Recall that $r_{\epsilon,i}(\Xbf)$ is defined as: \begin{equation}\label{r-epsi}
r_{\epsilon,i}(\Xbf):=\left\{\begin{array}{ll}
\frac{1}{2\epsilon}\|\gbf_i(\Xbf)\|^2,\quad \mbox{if}\quad \|\gbf_i(\Xbf)\|\le \epsilon,\\
\\
\|\gbf_i(\Xbf)\|-\frac{\epsilon}2,\quad \mbox{if}\quad \|\gbf_i(\Xbf)\|>\epsilon. 
\end{array}
\right.
\end{equation}

The proof of (\ref{third-a}) consists of three steps. In the first case, if $\|\gbf_i(\Xbf)\|\ge \delta$, it is easy to use the definition of $r_{\epsilon,i}(\Xbf)$ to see that both sides of (\ref{third-a}) are equal to
$\|\gbf_i(\Xbf)\|$. In the second case, we suppose $\epsilon<\|\gbf_i(\Xbf)\|<\delta$. Then
$$
  r_{\epsilon,i}(\Xbf)+\frac{\epsilon}{2}=\|\gbf_i(\Xbf)\|
  =\frac{\|\gbf_i(\Xbf)\|^2}{2\|\gbf_i(\Xbf)\|}+\frac{\|\gbf_i(\Xbf)\|}{2}
  \le \frac{\|\gbf_i(\Xbf)\|^2}{2\delta}+\frac{\delta}2=r_{\delta,i}(\Xbf)+\frac{\delta}2.
$$

Note that the last inequality holds due to the simple fact that the function 
$$h(t):=\frac{\|\gbf_i(\Xbf)\|}{2t}+\frac{t}2, $$
is nondecreasing for $t\ge \|\gbf_i(\Xbf)\|$. 

In case three: If $\delta>\epsilon\ge \|\gbf_i(\Xbf)\|$,
$$r_{\epsilon,i}(\Xbf)+\frac{\epsilon}2=\frac{1}{2\epsilon}\|\gbf_i(\Xbf)\|^2+\frac{\epsilon}2
\le \frac{1}{2\delta}\|\gbf_i(\Xbf)\|^2+\frac{\delta}2=r_{\delta,i}(\Xbf)+\frac{\delta}2
$$ where the inequality also holds due to the fact that
 the function 
$h(t)$ as defined above 
is a nondecreasing function for $\epsilon\ge \|\gbf_i(\Xbf)\|$. In each case, ${\bf C3}$ is satisfied.

{\bf Step three: Verification of {$\bf C4$}}:

Note that $\Xbf_{k_l+1}$ satisfies $\|\Phi_{\epsilon_{k_l}}(\Xbf_{k_l+1})\|\to 0$ as $l\to \infty$. For simplicity, we use $\{\Xbf_{j+1}\}$ instead of $\Xbf_{k_l+1}.$ $\epsilon_j$ is the smoothed factor used in the iteration to obtain $\Xbf_{j+1}$. Then, we have $\Xbf_{j+1}\to \tilde \Xbf$, $\epsilon_j\to 0$ and $\nabla_{1,2}\phi_{\epsilon_j}(\Xbf_{j+1})\to 0$. 

Here we claim that $\nabla_{1,2}\Phi_{\epsilon_j}(\Xbf_{j+1})$ converges to a Clarke point. Similar proof has been demonstrated in \cite{LDA} that 
\begin{align}\label{grad-phi}
\partial \Phi(\tilde \Xbf)=\{\sum_{i\in I_0}\nabla_{1,2} \gbf_i(\tilde \Xbf)^{\top}\wbf_i+\sum_{i\in I_1}\nabla _{1,2}\gbf_i(\tilde \Xbf)^{\top}\frac{\gbf_i(\tilde \Xbf)}{\|\gbf_i(\tilde \Xbf)\|}+\nabla_{1,2}f(\tilde \Xbf)  \bigg | \\
\wbf_{i}\in \mathbb C^d, \|\Pi(\wbf_i;\mathcal{C}(\nabla_{1,2} \gbf_i(\tilde \Xbf)))\|\le 1, \forall i\in I_0\} \nonumber
\end{align}
where $I_0=\{i\in [n] | \|\gbf_i(\tilde \Xbf)|=0\}$ and $I_1=[n]\setminus I_0$. $\Pi(\wbf;\mathcal{C}(A))$ is the projection of $\wbf$ onto $\mathcal{C}(A)$ which stands for the column space of $A$. Then for $J$ sufficiently large
$$\epsilon_j<\frac 12 \min\{\|\gbf_i(\tilde \Xbf)\|; i\in I_1\}\le \frac 12 \|\gbf_i(\tilde \Xbf)\|\le \|\gbf_i(\Xbf_{j+1})\|,\quad \forall j\ge J,\quad \forall i\in I_1 $$
We use the fact that $\min\{\| \gbf_i(\tilde \Xbf)\|; i\in I_1\}>0$ and $\epsilon_j\to 0$ to prove the first inequality and we use $\Xbf_{j+1}\to \tilde \Xbf$ and the continuity of $\gbf_i$ for all $i$ to prove the last inequality. Furthermore we denote 
$$\textbf{s}_{j,i}=\left\{\begin{array}{ll}\frac{\gbf_i(\Xbf_{j+1})}{\epsilon_j},\quad \mbox{if}\quad \|\gbf_i(\Xbf_{j+1})\|\le \epsilon_j,\\
\\
\frac{\gbf_i(\Xbf_{j+1})}{\|\gbf_i(\Xbf_{j+1})\|},\quad \mbox{if}\quad \|\gbf_i(\Xbf_{j+1})\|>\epsilon_j.
\end{array}
\right. $$
Then, we have 
$$\nabla_{1,2} \Phi_{\epsilon_j}(\Xbf_{j+1})=\sum_{i\in I_0}\nabla_{1,2}\gbf_i(\Xbf_{j+1})^{\top}\textbf{s}_{j,i}+\sum_{i\in I_1}\nabla_{1,2} \gbf_i(\Xbf_{j+1})^{\top}\frac{\gbf_i(\Xbf_{j+1})}{\|\gbf_i(\Xbf_{j+1})\|}+\nabla_{1,2}f(\Xbf_{j+1}).$$
As $j\to \infty,$ obviously $\nabla_{1,2} \Phi_{\epsilon_j}(\Xbf_{j+1})$ converges to the three terms in (\ref{grad-phi}). Since $\nabla_{1,2} \Phi_{\epsilon_j}(\Xbf_{j+1})\to 0$ and $\partial \Phi(\tilde \Xbf)$ is closed, $\Xbf_{j+1}$ converges to a Clarke stationary point. We have validated the four conditions required for the proposed algorithm to be convergent in regards of the application of joint reconstruction of T1 and T2 MRI images. Therefore, by Lemma 1 and Theorem 1, this joint reconstruction will provide a convergent solution.

\subsection{Experiment settings}
The dataset we used in all experiments are from {\it Multi-modal Brain Tumor Segmentation Challenge 2018}.\cite{menze2014multimodal} Both the training set and the validation set in Multi-modal Brain Tumor Segmentation Challenge 2018 contain four modalities (T1,T2, FLAIR and T$1_{c}$). The training set is scanned from 285 patients and the validation set is obtained from 66 patients. Each image has a volume size of $240\times 240\times 155.$ Each modality has two types of gliomas: 75 volumes of low-grade gliomas(LGG) and 210 volumes of high-grade gliomas(HGG). 

In our experiments, we used HGG MRI images from two modalities (T1 and T2). For each modality, we chose images scanned from 50 patients as our training set and we randomly chose images scanned from 6 patients as our testing set. We cropped the 2D image size to be $160\times 180$ in the center region, resulting in a overall number of 500 images as our training set and 60 images as our testing set. 

To obtain the k-space images, We used Matlab to apply 2D fast Fourier transform to every ground truth image in the training set and the testing set. Then, we shifted the zero-frequency components to the central area for each image and we applied 10\%, 20\% radial mask to every image. After having all the operations on Matlab, we had the k-space under-sampled images, denoted as $\textbf{f}_{1}$ and $\textbf{f}_{2}$.  

The Initialization network was implemented in Python using the TensorFlow framework and the LPAM-net was implemented in Python using the PyTorch framework. Our experiments were run on a Linux server with an NVIDIA A100
Tensor Core GPU available on HiPerGator. The batch size was set to be 1 in all experiments due to the consideration of the GPU memory, data volume and the operation speed.

Our reconstruction results are assessed using four evaluation matrices: peak-to-noise ratio (PSNR), structural similarity (SSIM), normalized mean squared error (NMSE) and root mean square error (RMSE). PSNR, SSIM, NMSE and RMSE can be computed using the following equations.

Suppose we have a reconstructed image $\Xbf$ and a ground truth image $\textbf{Y}$, each with a dimension of $ m\times n.$ We define the mean square error and peak-to-noise ratio as 
\begin{equation}\label{mse}
MSE \hspace{0.5mm}(\Xbf,\textbf{Y})=\dfrac{1}{mn}\sum_{i=0}^{m-1}\sum_{j=0}^{n-1}[\Xbf(i,j)-\textbf{Y}(i,j)]^{2}=\dfrac{1}{mn}||\Xbf-\textbf{Y}||_{2}^{2}
\end{equation}

$$PSNR \hspace{0.5mm}(\Xbf,\textbf{Y})=10\cdot log_{10}\left(\dfrac{MAX(\textbf{Y})}{MSE(\Xbf,\textbf{Y})}\right)$$
where $MAX(\textbf{Y})$ is the maximum possible pixel value of the ground truth image. PSNR measures the quality of a compressed image in comparison to its ground truth image. A higher PSNR indicates a better quality of the compressed or reconstructed image.

$$SSIM\hspace{0.5mm}(\Xbf,\textbf{Y})=\dfrac{(2\mu_{X}\mu_{Y}+c_{1})(2\sigma_{XY}+c_{2})}{(\mu_{X}^{2}+\mu_{Y}^{2}+c_{1})(\sigma_{X}^{2}+\sigma_{Y}^{2}+c_{2})}$$
where $\mu_{X},\mu_{Y}$ is the pixel mean of $X$ and $Y,$ respectively; $\sigma_{X}^{2},\sigma_{Y}^{2}$ are the variance of $X, Y$ respectively. $c_{1}=(k_{1}L)^{2},c_{2}=(k_{2}L)^{2}$ are two variables to avoid zero denominator and $L$ is the dynamic range of the pixel values. The SSIM, ranging from 0 to 1, is used to measure the similarity between two images. A higher SSIM value indicates a greater similarity between the two images.

$$NMSE\hspace{0.5mm}(\Xbf,\textbf{Y})=\dfrac{||\Xbf-\textbf{Y}||_{2}^{2}}{||\textbf{Y}||_{2}^{2}}$$
NMSE measures the mean relative error. 

$$RMSE\hspace{0.7mm}(\Xbf,\textbf{Y})=\sqrt{MSE(\Xbf,\textbf{Y})}=\sqrt{\dfrac{1}{mn}||\Xbf-\textbf{Y}||_{2}^{2}}$$
RMSE is utilized to measure the difference between the reconstructed image and the ground truth. 

\subsection{The architecture of the reconstruction process }
The complete reconstruction procedure consists of the Initialization network and the LPAM-net. The architecture of the multi-phase LPAM-net follows the LPAM algorithm exactly in the way that each phase of the network corresponds to each itration of the algorithm. The details are given in the following subsections.

\subsubsection{Initialization network}

To have a better input for LPAM-net, we construct an initialization network, the major component of which is the initialization blocks.
It contains a 4-layer convolutional neural network(CNN) following the residual structure, as outlined in \cite{he2016deep}. Subsequently, an inverse Fourier transform is applied to the output of the initialization block. 
The initialization block performs interpolation on the inputs $\textbf{f}_{1}$ and $\textbf{f}_{2}$, which helps to fill up the missing components of the under-sampled images in the k-space domain. We learned the weights of the kernels of the Initialization network by minimizing the following loss function $L=Loss1+Loss2$ using Adam Optimizer implemented using Tensorflow:
  $$ Loss1=\dfrac{1}{mn}||\xbf_{1}^{0}-\xbf_{1}^{truth}||_{2}^{2}$$
  $$  Loss2=\dfrac{1}{mn}||\xbf_{2}^{0}-\xbf_{2}^{truth}||_{2}^{2},$$
where $m\times n$ is the dimension of the image.


A learning rate of $10^{-3}$ was set, complemented by a decay rate of 0.95 for every 100 steps. The Xavier initialization was utilized to initialize the kernels. Table \ref{tab:ini} displays the quantitative assessment of the Initialization network for the 60 images in the testing set. This includes the average and the standard deviation for the outputs {$\xbf_{1}^{0}, \xbf_{2}^{0}$} under noiseless condition, corresponding to radial masks with an under-sampling ratios of 10\% and 20\%.

\begin{table}[h]
\centering
\addtolength{\tabcolsep}{-1pt}
\caption{The evaluation results (mean $\pm$ standard deviation) for noiseless $\xbf_{1}^{0}$ and $\xbf_{2}^{0}$ after the initialization network  }

\begin{tabular}{lcccc}
\toprule
\textbf{Image} & \textbf{PSNR} & \textbf{SSIM} & \textbf{NMSE} \\
\midrule 
10\% T1 & 23.146 $\pm$ 0.66 & 0.617 $\pm$ 0.034 & 0.0196 $\pm$ 0.0035 \\
10\% T2 & 23.65 $\pm$ 1.64 & 0.619 $\pm$ 0.04 & 0.063 $\pm$ 0.026 \\
20\% T1 & 26.517 $\pm$ 0.518 & 0.735 $\pm$ 0.031 & 0.009 $\pm$ 0.0014  \\
20\% T2& 27.032 $\pm$ 1.633 & 0.735 $\pm$ 0.031 & 0.029 $\pm$ 0.012 \\

\bottomrule
\toprule
\end{tabular}

\label{tab:ini}
\end{table}

\subsubsection{LPAM-net for Joint Reconstruction }
\begin{figure}[H]
    \centering
    \includegraphics[width=1\textwidth]{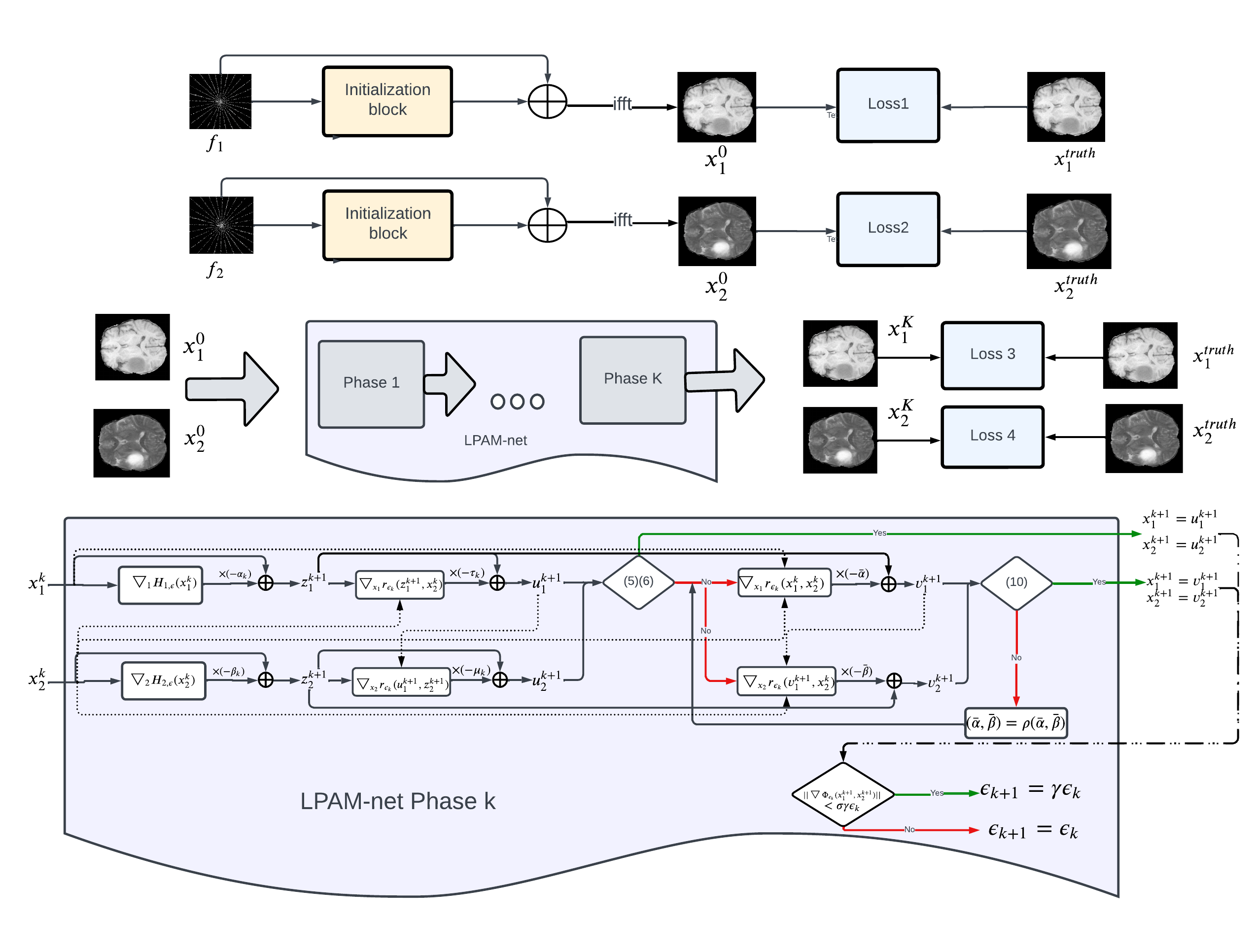}
    \caption{The architecture of the propsoed Initialization network and the LPAM-net. Top: The proposed Initialization network; Middle: The architecture of the LPAM-net for joint reconstruction of T1 and T2 images; Bottom: The detailed illustration of the $k^{th}$--phase of the LPAM-net.}
    \label{fig:joint}
\end{figure}

The architecture of the proposed LPAM-net  is shown in the flowchart \ref{fig:joint}, which follows the LPAM algorithm exactly. The inputs for the LPAM-net are the outputs of the Initialization network, denoted as $\xbf_{1}^{0}$ and $\xbf_{2}^{0}$. The outputs of the final phase of the LPAM-net are the resulting reconstructed images.

Throughout our experiments, the regularization we used is the $\ell_{2,1}$-norm of the joint feature extractor, which is described in the Section 7.1. The formula for the regularization term is $r(\xbf_{1},\xbf_{2})=||\gbf_{\theta}(\xbf_{1},\xbf_{2})||_{2,1},$ and the expressions of $\nabla_{1,2} r_{\epsilon}(\xbf_1,\xbf_2)$ and $r_{\epsilon}(\xbf_1,\xbf_2)$ are given in \ref{deltar} and \ref{r-epsi}, respectively. We denote the feature extractor $\gbf_{\theta}(\xbf_{1},\xbf_{2})$ as $\gbf(\xbf_{1},\xbf_{2})$ or $\gbf$  for simplicity in the following text. Regarding the joint feature extractor $\gbf$, we parameterize it as \ref{g} with component-wise activation $\sigma$ as \ref{sigma}. The prefixed parameter $\delta$ in the activation function $\sigma$ was set to be 0.01. The default configuration of the feature extraction operator $\gbf$ was set as follows: it consisted of $l=4$ complex-kernel convolution layers.
Every kernel in each layer has
a dimension of 3$\times$3$\times$32 (except for the first layer), where 32 is the depth of the kernel. 
In the first layer, the kernel size was 3$\times$3$\times$2, with a depth of 2. For each layer, we had 32 kernels. We set the stride to be 1 and used zero-padding to preserve image size. 

For joint reconstruction using the LPAM-net, we need to learn the step sizes $\alpha_{k},\tau_{k},\beta_{k}$ and $\gamma_{k},$ the threshold $\epsilon_{k},$ and the weights of the kernels of the feature extractor $\gbf$. We consolidated all the learned parameters into a unified symbol denoted as $\theta$. The total number of the phases k is set to be $15$.
To enhance the stability and precision of the parameters before increasing the iteration count, 
we initially trained LPAM-net with phase number K=3 for 100 epochs. After completing the training of the phase K=3, an incremental strategy involving adding 2 more phases each time was employed and for K$>$3, we conducted training of 30 epochs during each iteration. Throughout the training, the batch size was set to be $1$ and parameters were updated after the training for every batch.

We updated the parameters $\theta$ by minimizing the loss function $L'$ using Adam Optimizer implemented by Pytorch with $\beta_{1}=0.9, \beta_{2}=0.999$. The learning rate was set at $10^{-4}$, and the initial step sizes were established as $\alpha_{0}=\beta_{0}=0.5$. For $\tau_{0}$ and $\mu_{0},$ we have three  initial step sizes for different phases. For phase 1 to phase 3, $\tau_{0}=\mu_{0}=2$; for phase 4 to phase 12, $\tau_{0}=\mu_{0}=1;$ for phase 13 to phase 15, $\tau_{0}=\mu_{0}=0.1.$ The Xavier initialization was applied to initialize the kernels. The outputs for the algorithm were denoted as $\xbf_{1}^{k}$ and $\xbf_{2}^{k}$ and the ground truth were represented by $\xbf_{1}^{truth}$ and $\xbf_{2}^{truth}.$ To calculate the loss, we converted the $m\times n$ image into a vector of length $mn$. Suppose the total pixel number of an image is $mn,$ the loss function we minimized is expressed as $L'=Loss3+Loss4$ where $Loss3$ and $Loss4$ are defined as follows respectively:\\
$$ Loss3=\dfrac{1}{mn}||\xbf_{1}^{k}-\xbf_{1}^{truth}||_{2}^{2}+0.1\times(1-SSIM(\xbf_{1}^{k},\xbf_{1}^{truth}))$$
$$  Loss4=\dfrac{1}{mn}||\xbf_{2}^{k}-\xbf_{2}^{truth}||_{2}^{2}+0.1\times(1-SSIM(\xbf_{2}^{k},\xbf_{2}^{truth}))$$

\subsection{Convergence behavior of the LPAM-net}

\begin{figure}[H]
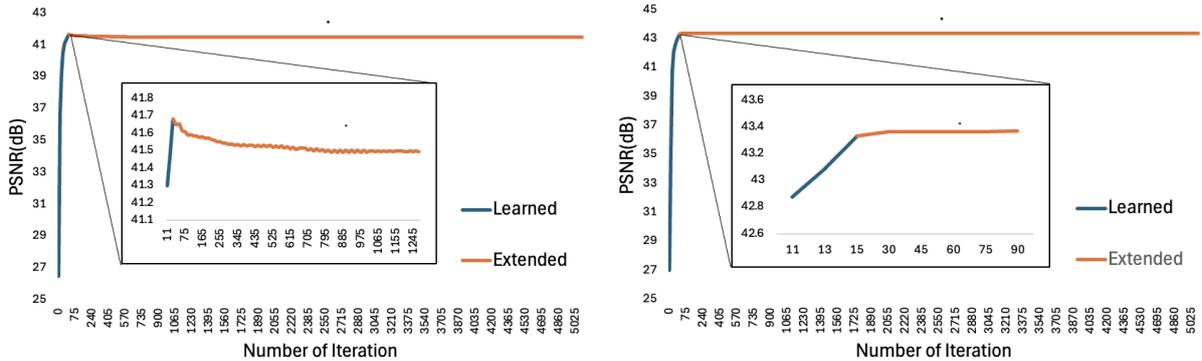

\includegraphics[width=0.5\textwidth]{t150000.pdf}
\includegraphics[width=0.5\textwidth]{t250000.pdf}
\caption{From Left to Right: The PSNR values for the reconstruction results of the T1 and T2 image with a $20\%$ under-sampling ratio at various phases.  In both figures, the LPAM-net  learns the network parameters in the first 15 phases.  Starting from $\Xbf_{16}$ the reconstruction results are obtained from the same model \eqref{modell} with fixed $\gbf$ learned from the 15-phase LPAM-net. The PSNR values for the reconstruction results  of the first 15 phases are plotted in blue, and since the $16^{th}$ phase, that are plotted in orange.}
\label{yuo}
\end{figure}

In the previous subsection, we proved that there is at least a subsequence of the iterates generated by the LPAM algorithm converges to a Clarke stationary point. It is expected that LPAM would perform stably even beyond the trained phases. 
To demonstrate its stability empirically,
we set the parameters as follows: $\epsilon_{0}=0.01,$ $\gamma=0.9,$ $\sigma=60000$ and the number of phases $K$ = 15. We used the $k$-space data of a $20\% $ under-sampling ratio for the joint reconstruction of MR T1 and T2 images. The $\ell_{2}$-norm of the gradient of the objective function value  ($\|\nabla_{1,2}\Phi_{\varepsilon_k}(\xbf_1^{k}, \xbf_2^{k})\|$) decreases as the number of phase increases. According to the $\epsilon$ reduction criteria specified in the line 19 in the algorithm, this leads to a reduction in the smoothing factor from $0.01$ to $0.0043$ by phase $15$.
During the 15 phases, the smoothing factor steadily decreases as the number of phases increases, indicating that the smoothed objective function is progressively converging to the original problem.

\begin{figure}[H]

    \centering
    \includegraphics[width=0.6\textwidth]{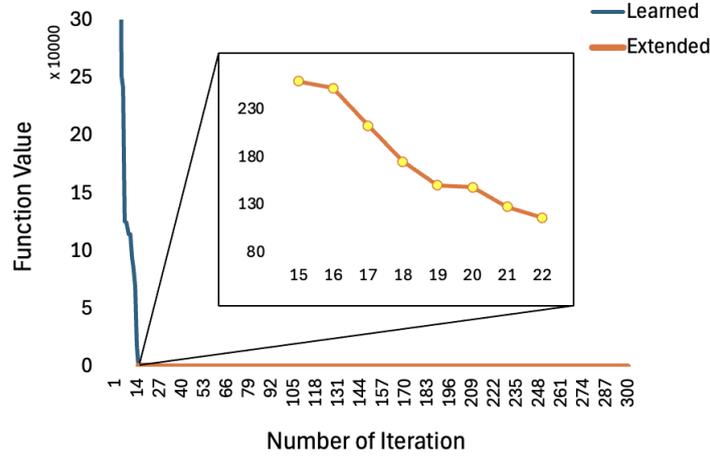}
    
    \caption{Objective function value ($\Phi(\xbf_1^{k}, \xbf_2^{k})$) of images reconstruction with a 20\% under-sampling ratio across various phase number K. The results of the first 15 phases are plotted in blue and since the $16^{th}$ phase, that are plotted in orange.}
    \label{fig:phii}
\end{figure}

\begin{figure}[H]
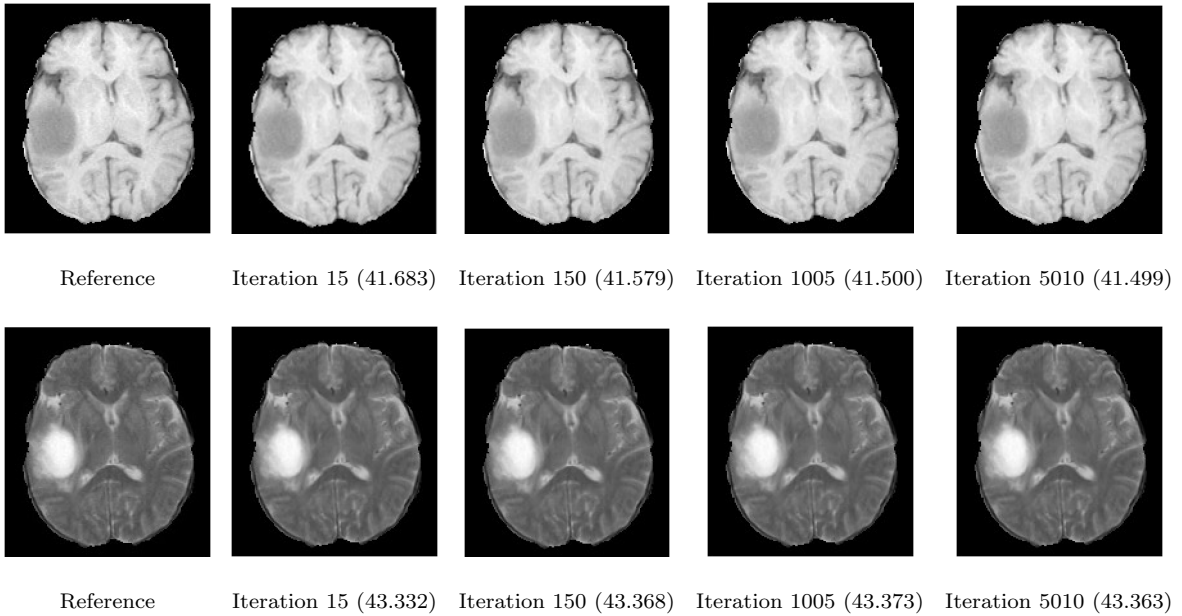


  \centering
  \setlength{\tabcolsep}{4pt}
  \begin{tabular}{c c c c c}
    \includegraphics[width=0.17\textwidth]{T1_TRUE_.pdf} &
    \includegraphics[width=0.17\textwidth]{T1_15_.pdf} &
    \includegraphics[width=0.17\textwidth]{T1_150_.pdf}&
    \includegraphics[width=0.17\textwidth]{T1_1005_.pdf} &
    \includegraphics[width=0.17\textwidth]{T1_5010_.pdf} 
    
    \\
    [\abovecaptionskip] \small Reference &
    \small  Iteration 15 (41.683)&
    \small Iteration 150 (41.579) &
    \small Iteration 1005 (41.500) &
    \small Iteration 5010 (41.499)
  \end{tabular}

  \vspace{\floatsep}

  \begin{tabular}{c c c c c}
    \includegraphics[width=0.17\textwidth]{T2_TRUE_.pdf} &
    \includegraphics[width=0.17\textwidth]{T2_15_.pdf} &
    \includegraphics[width=0.17\textwidth]{T2_150_.pdf} &
    \includegraphics[width=0.17\textwidth]{T2_1005_.pdf} &
    \includegraphics[width=0.17\textwidth]{T2_5051.pdf} 
    \\
    [\abovecaptionskip] \small Reference &
    \small Iteration 15 (43.332) &
    \small Iteration 150 (43.368) &
    \small Iteration 1005 (43.373)&
    \small Iteration 5010 (43.363)
  \end{tabular}

  \caption{
Top row: Representation brain T1 MRI images reconstructed by the LPAM-net with a 20\% under-sampling ratio in the radial mask after 15, 150, 1005 and 5010 iterations; Bottom row: Representation brain T2 MRI images reconstructed by the LPAM-net with a 20\% under-sampling ratio in radial mask after 15, 150, 1005 and 5010 iterations. PSNRs (dB) are shown in the parentheses. }
\label{fig:iteration}
\end{figure}

Furthermore, 
we extend the algorithm beyond 15 phases and observe its performance on the values of the objective function and PSNR for the reconstruction. To this end,
after $15$ phases, we set $\epsilon_{tol}=0$ and ran the algorithm for up to 5010 iterations with fixed $g$ that is learned at the $15^{th}$ phase. Figure \ref{yuo} and figure \ref{fig:phii} illustrate the changes in the values of PSNR and objective function as the number of iterations increases, respectively. Figure \ref{fig:phii} is drawn based on the objective function with 0.0093 as the coefficient of the regularization part. This coefficient does not affect the convergence based on Lemma 1. As we can tell, LPAM-net reduces the objective function value steadily in the first 15 phases. After that, the value of the objective function continue to decrease, but not as fast as the first 15 phases. 
The PSNR for both T1 and T2 images after the $15^{th}$ phase maintain stable with slightly decrease in the first a few iterations. The figure \ref{fig:iteration} visually  shows the reconstructed images obtained after 15, 150, 1005, 5010 iterations. We can see that the reconstructed images are almost identical and free of visual artifacts, with very similar PSNR for both T1 and T2 contracts. These evidence shows that LPAM-net maintains stable performance over extended iterations.


\subsection{Results}
In this subsection, we provide numerical results to showcase the enhanced image quality achieved by the proposed LPAM-net and the variational model \ref{modell}. First, in order to show the advantage of using the joint features as the regularization term in the model, we make comparison between the LPAM-net and the Individual-modality Reconstruction Network. Next, we compare the LPAM-net with the network induced by BCD algorithm, which is equivalent to run only the lines 10-17 and skip  the lines 3-9 in LPAM algorithm. Finally, we present the comparison of the LPAM-net with four state-of-art methods for joint reconstruction of MRI T1 and T1 images.

\subsubsection{Comparison with Individual-modality Reconstruction Network}
This experiment is designed to verify that using the joint features in the regularization term enhances the accuracy and efficiency of the reconstructions of the MRI T1 and T2 images. Since the LPAM-net takes advantage of the common features from both T1 and T2 contracts, the proposed reconstructed images should preserve more details. To validate this advantage, we designed an Individual-modality Reconstruction Network which extracts the individual features of the T1 and T2 contracts separately. The Individual-modality Reconstruction Network precisely adheres the algorithm in \cite{LDA} for minimizing the following problem:

$$\min_{\xbf_{i}} \|\textbf{P}\textbf{F}\xbf_i-\textbf{f}_i\|^2+\|\gbf_{\theta_{i}}(\xbf_i)\|_{2,1} $$
where $\xbf_i\in\mathbb R^n$, $\textbf{f}_i$ is the under-sampled k-space images and $i=1,2.$ $\gbf_{\theta_{i}}$ is the individual feature extractor from $\mathbb C^n$ to $\mathbb C^{n\times d}$ and it is designed to be a 4-layer CNN. Unlike the joint reconstruction, in which we learn a joint feature extraction operator $\gbf_{\theta}$, in this individual reconstruction, we learn two feature extraction operators $\gbf_{\theta_1}$ and $\gbf_{\theta_2}$ 
separately. We learned the weights of the kernels in the CNNs by minimizing the loss function using Adam Optimizer implemented by Pytorch. Suppose the total pixel number of an image is $mn,$ the loss function we minimized is expressed as follows:
$$ Loss=\dfrac{1}{mn}||\xbf_{i}^{k}-\xbf_{i}^{truth}||_{2}^{2}+0.1\times(1-SSIM(\xbf_{i}^{k},\xbf_{i}^{truth}))$$
where $\xbf_{i}^{k}$ is the output and $\xbf_{i}^{truth}$ is the ground truth image.



\begin{figure}[H]
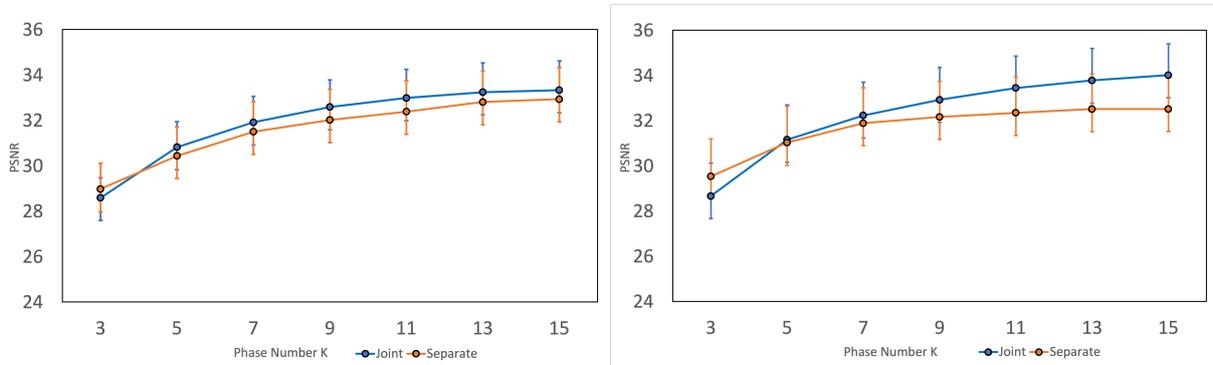

\includegraphics[width=0.5\textwidth]{T1ff.pdf}
\includegraphics[width=0.5\textwidth]{T2ff.pdf}
\caption{ PSNR for the reconstruction result of images with a 10\% under-sampling ratio  across various phase number $K$. Left: T1 images. Right: T2 images. 
In both figures, the blue line  represents results for joint reconstruction via LPAM-net and the orange line  represents the results from the Individual-modality Reconstruction Network. The vertical line represents the standard deviation of PSNR.}
\label{psnr-phase}
\end{figure}

We compare T1 and T2 image reconstruction results with  $10\%$ under-sampling ratio $k$-space noiseless data resulting from LPAM-net and the Individual-modality reconstruction Network, respectively. Figure \ref{psnr-phase}  shows the PSNR value at the phase number $K=3,5,...,15$ of both deep neural networks. Both networks were trained using Adam Optimizer implemented in Pytorch with $\beta_{1}=0.9, \beta_{2}=0.999.$ The learning rate was set at $10^{-4}$ and the number of phase is set at 15.
The average PSNR and its standard deviation are depicted in figure \ref{psnr-phase} for each phase. For phase 3, Individual-modality Reconstruction Network outperforms, but as the phase number increases, the LPAM-net produces better reconstructions. At phase 15, the average PSNR of the LPAM-net improves 0.40 dB for the T1 images and improves 1.49 dB for the T2 images, respectively comparing to the Individual-modality Reconstruction Network. In addition, the standard deviation decreases 0.10 dB for the T1 images and decreases 0.21 dB for the T2 images.

Although PSNR is a crucial measure for comparing the two networks, it is hard for us to visually discern the impact that the LPAM-net has in the reconstruction process. In order to visualize the differences, we added Gaussian white noise to the complex-valued under-sampled k-space images. The real noise and the imaginary noise are both normally distributed with the same mean and the same standard deviation, but they maintain independence from each other. Then, we acquired three groups of T1, T2 images: one set is noiseless; one set has complex-valued noise with a mean of zero and a standard deviation of 3 on both the real part and imaginary part, and the third set has complex-valued noise with a mean of zero and a standard deviation of 7 on both the real part and imaginary part. Next, we trained and tested these three groups of T1,T2 images using both the LPAM-net and the Individual-modality Reconstruction Network. The resulting reconstructed T1 images from the three groups are presented in figure \ref{fig:t1oo} with the ground truth reference. Figure \ref{fig:t1dd} illustrates the pointwise absolute error for each group of the images compared to the ground truth reference. Likewise, the reconstructed T2 images and the associated absolute errors are shown in figure \ref{fig:t2qq} and figure \ref{fig:t2dd}. 

Figure \ref{fig:t1oo} and figure \ref{fig:t2qq} reveal that, for under-sampled images with any levels of Gaussian white noise, the edges of the tissues 
in the images reconstructed by the LPAM-net are 
sharper than the edges of the tissue 
reconstructed by the Individual-modality Reconstruction Network. Furthermore, 
the absolute errors are less in the images reconstructed from the LPAM-net in figure \ref{fig:t1dd} and figure \ref{fig:t2dd}. Table \ref{table11} and table \ref{table22} provide a comprehensive summary for the quantitative 
comparison between the LPAM-net and the Individual-modality Reconstruction Network across various levels of noise. The PSNR, SSIM, NMSE and RMSE are provided for images with each level of noise. We can conclude that the LPAM-net generates more accurate images than the Individual-modality Reconstruction Network regardless of the noise level in the under-sampled images. The results show the advantage of using the joint feature operator as the regularization term in the reconstruction process. Additionally, using joint feature operator demonstrates greater parameter efficiency.

\begin{table}[H]
\centering
\addtolength{\tabcolsep}{-1pt}
\caption{The results (mean $\pm$ standard deviation) of T1 image reconstruction using LPAM-net and Individual-modality Reconstruction Network, with the radial mask at a 10\% under-sampling ratio, are presented for different variances of Gaussian white noises added to the k-space data.}
\label{table11}
\begin{tabular}{lcccccc}
\toprule
   \textbf{Method}& \textbf{SD of noise}&
   \textbf{PSNR} & \textbf{SSIM} & \textbf{NMSE} &
   \textbf{RMSE}&\textbf{\#Par}\\
\midrule 

Individual-modality & 7 & 29.1356 $\pm$ 0.9087 & 0.8688 $\pm$ 0.0108 & 0.0051 $\pm$ 0.0016 & 0.0351 $\pm$ 0.0036& 55903\\
LPAM-net & 7 & 29.9429 $\pm$ 0.7945 & 0.8852 $\pm$ 0.0097 & 0.0042 $\pm$ 0.0011 & 0.032 $\pm$ 0.0028& 56510\\
Individual-modality & 3 & 31.2754 $\pm$ 1.1333 & 0.9044 $\pm$ 0.0098 & 0.0031 $\pm$ 0.0011 & 0.0275 $\pm$ 0.0034& 55903\\
LPAM-net & 3 & 31.994 $\pm$ 1.01 & 0.9175 $\pm$ 0.0091 & 0.0026 $\pm$ 0.0008 & 0.0253 $\pm$ 0.0028 & 56510\\
Individual-modality & none & 33.0252 $\pm$ 1.3921 & 0.9251 $\pm$ 0.0094 & 0.0021 $\pm$ 0.0008 &0.0226 $\pm$ 0.0034 & 55903\\
LPAM-net & none & 33.3269 $\pm$ 1.2853 & 0.927 $\pm$ 0.0091 & 0.002 $\pm$ 0.0007 & 0.0218 $\pm$ 0.003 & 56510\\

\bottomrule
\toprule
\end{tabular}
\end{table}

\begin{table}[htpb]
\centering
\addtolength{\tabcolsep}{-1pt}
\caption{The results (mean $\pm$ standard deviation) of T2 image reconstruction using LPAM-net and Individual-modality Reconstruction Network, respectively, with the radial mask at a 10\% under-sampling ratio, are presented for different variances of Gaussian white noises added to the k-space data.}
\label{table22}

\begin{tabular}{lcccccc}

\toprule
   \textbf{Method}& \textbf{SD of noise}&
   \textbf{PSNR} & \textbf{SSIM} & \textbf{NMSE} &
   \textbf{RMSE}&\textbf{\#Par} \\

\midrule 

Individual-modality & 7 & 28.6758 $\pm$ 1.0697 & 0.8632 $\pm$ 0.016 & 0.0193 $\pm$ 0.0067 & 0.0371 $\pm$ 0.0046 & 55903\\
LPAM-net & 7 & 29.8219 $\pm$ 1.107 & 0.88 $\pm$ 0.0178 & 0.015 $\pm$ 0.0055 & 0.0325 $\pm$ 0.0041& 56510\\
Individual-modality & 3 & 30.972 $\pm$ 1.2514 & 0.9049 $\pm$ 0.0149 & 0.0114 $\pm$ 0.0040 & 0.0286 $\pm$ 0.0041& 55903\\
LPAM-net & 3 & 32.1759 $\pm$ 1.0952 & 0.9209 $\pm$ 0.0147 & 0.0087 $\pm$ 0.0031 & 0.0248 $\pm$ 0.0031 & 56510\\
Individual-modality & none & 32.621 $\pm$ 1.6011 & 0.9209 $\pm$ 0.0151 & 0.0078 $\pm$ 0.0028 &0.0238 $\pm$ 0.0044  & 55903\\
LPAM-net & none & 34.011 $\pm$ 1.3859 & 0.9392 $\pm$ 0.0132 & 0.0057 $\pm$ 0.0021 & 0.0202 $\pm$ 0.0032 & 56510\\
\bottomrule
\toprule
\end{tabular}
\end{table}

\begin{figure}[htpb]
  \centering
  \setlength{\tabcolsep}{4pt}
  \begin{tabular}{c c c c}
    \includegraphics[width=0.18\textwidth,angle=-180,clip]{T17sepoo.pdf} &
    \includegraphics[width=0.18\textwidth,angle=-180,clip]{T13sepoo.pdf} &
    \includegraphics[width=0.18\textwidth,angle=-180,clip]{T1nonsepoo.pdf} &
    \includegraphics[width=0.18\textwidth,angle=-180,clip]{T1trueoo.pdf} 
    \\
    [\abovecaptionskip] \small Individual\_7 (29.233) &
    \small  Individual\_3 (31.373)&
    \small Individual\_non (32.933) &
    \small Reference
  \end{tabular}

  \vspace{\floatsep}
  
  \setlength{\tabcolsep}{4pt}
  \begin{tabular}{c c c c}
    \includegraphics[width=0.18\textwidth,angle=-180,clip]{T17jointoo.pdf} &
    \includegraphics[width=0.18\textwidth,angle=-180,clip]{T13jointoo.pdf} &
    \includegraphics[width=0.18\textwidth,angle=-180,clip]{T1nonjointoo.pdf} &
    \includegraphics[width=0.18\textwidth,angle=-180,clip]{T1trueoo.pdf} 
    \\
    [\abovecaptionskip] \small LPAM-net\_7 (29.943) &
    \small LPAM-net\_3 (31.994) &
    \small LPAM-net\_non (33.327) &
    \small Reference
  \end{tabular}

  \caption{Reconstruction results of T1 brain images
by the LPAM-net and the Individual-modality Reconstruction Network employing the radial mask with a 10\% under-sampling ratio. Top to bottom rows: Reconstruction results obtained
by the Individual-modality Reconstruction Network and
the LPAM-net, respectively. The first to the third columns: Images reconstructed from the k-spaces, where a complex noise with a standard deviation of 7, 3, and  0 (noiseless) is added, respectively.
The last column: Ground truth reference. PSNRs (dB) are shown in the parentheses. The regions of interest are magnified in red boxes for better visualization.}
\label{fig:t1oo}
\end{figure}

\begin{figure}[htpb]
  \centering
  \setlength{\tabcolsep}{4pt}
  \begin{tabular}{c c c }
    \includegraphics[width=0.18\textwidth,angle=-180,clip]{T17sepdd.pdf} &
    \includegraphics[width=0.18\textwidth,angle=-180,clip]{T13sepdd.pdf} &
    \includegraphics[width=0.18\textwidth,angle=-180,clip]{T1nonsepdd.pdf} 
    \\
    [\abovecaptionskip] \small Individual\_7 (29.233) &
    \small  Individual\_3 (31.373)&
    \small Individual\_non (32.933) 
  \end{tabular}

  \vspace{\floatsep}

  \begin{tabular}{c c c }
    \includegraphics[width=0.18\textwidth,angle=-180,clip]{T17jointdd.pdf} &
    \includegraphics[width=0.18\textwidth,angle=-180,clip]{T13jointdd.pdf} &
    \includegraphics[width=0.18\textwidth,angle=-180,clip]{T1nonjointdd.pdf} 
    \\
    [\abovecaptionskip] \small LPAM-net\_7 (29.943) &
    \small LPAM-net\_3 (31.994) &
    \small LPAM-net\_non (33.327) 
  \end{tabular}

  \caption{
The corresponding pointwise absolute errors between the reconstructed brain T1 images in Figure \ref{fig:t1oo} and the ground truth reference. The errors are scaled up by a factor of 4 for better visualization and the bright parts indicate large value. 
}
\label{fig:t1dd}
\end{figure}

\begin{figure}[htpb]
  \centering
  \setlength{\tabcolsep}{4pt}
  \begin{tabular}{c c c c }
    \includegraphics[width=0.18\textwidth,angle=-180,clip]{T27sepoo.pdf} &
    \includegraphics[width=0.18\textwidth,angle=-180,clip]{T23sepoo.pdf} &
    \includegraphics[width=0.18\textwidth,angle=-180,clip]{T2nonsepoo.pdf} &
    \includegraphics[width=0.18\textwidth,angle=-180,clip]{T2trueoo.pdf} 
    \\
    [\abovecaptionskip] \small Individual\_7 (28.989) &
    \small  Individual\_3 (31.052)&
    \small Individual\_non (32.521) &
    \small Reference
  \end{tabular}

  \vspace{\floatsep}

  \begin{tabular}{c c c c }
    \includegraphics[width=0.18\textwidth,angle=-180,clip]{T27jointoo.pdf} &
    \includegraphics[width=0.18\textwidth,angle=-180,clip]{T23jointoo.pdf} &
    \includegraphics[width=0.18\textwidth,angle=-180,clip]{T2nonjointoo.pdf} &
    \includegraphics[width=0.18\textwidth,angle=-180,clip]{T2trueoo.pdf} 
    \\
    [\abovecaptionskip] \small LPAM-net\_7 (29.822) &
    \small LPAM-net\_3 (32.176) &
    \small LPAM-net\_non (34.011) &
    \small Reference
  \end{tabular}

   \caption{Reconstruction results of T2 brain images
by the LPAM-net and the Individual-modality Reconstruction Network employing the radial mask with a 10\% under-sampling ratio. Top to bottom rows: Reconstruction results obtained
by the Individual-modality Reconstruction Network and
the LPAM-net, respectively. The first to the third columns: Images reconstructed from the k-spaces, where a complex noise with a standard deviation of 7, 3, and  0 (noiseless) is added, respectively.
The last column: Ground truth reference. PSNRs (dB) are shown in the parentheses. The regions of interest are magnified in red boxes for better visualization.}

\label{fig:t2qq}
\end{figure}

\begin{figure}[htpb]
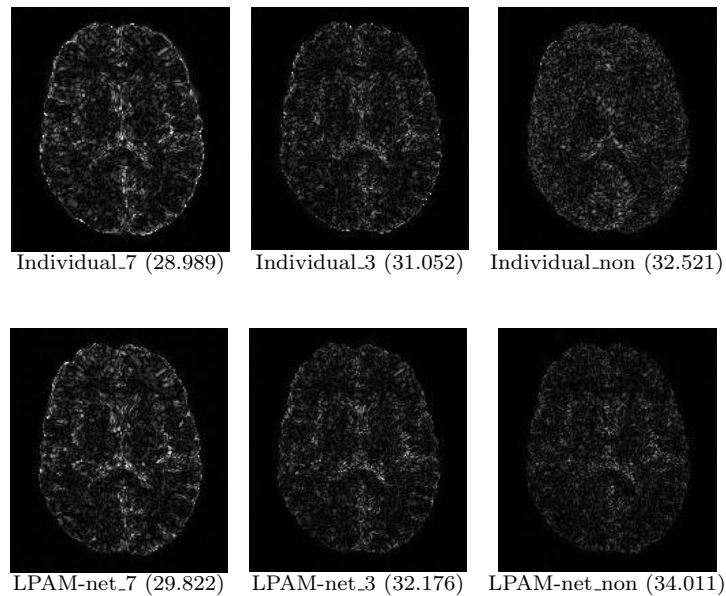

  \centering
  \setlength{\tabcolsep}{4pt}
  \begin{tabular}{c c c }
    \includegraphics[width=0.18\textwidth,angle=-180,clip]{T27sepdd.pdf} &
    \includegraphics[width=0.18\textwidth,angle=-180,clip]{T23sepdd.pdf} &
    \includegraphics[width=0.18\textwidth,angle=-180,clip]{T2nonsepdd.pdf} 
    \\
    [\abovecaptionskip] \small Individual\_7 (28.989) &
    \small  Individual\_3 (31.052)&
    \small Individual\_non (32.521) 
  \end{tabular}

  \vspace{\floatsep}

  \begin{tabular}{c c c }
    \includegraphics[width=0.18\textwidth,angle=-180,clip]{T27jointdd.pdf} &
    \includegraphics[width=0.18\textwidth,angle=-180,clip]{T23jointdd.pdf} &
    \includegraphics[width=0.18\textwidth,angle=-180,clip]{T2nonjointdd.pdf} 
    \\
    [\abovecaptionskip] \small LPAM-net\_7 (29.822) &
    \small LPAM-net\_3 (32.176) &
    \small LPAM-net\_non (34.011) 
  \end{tabular}

  \caption{
The corresponding pointwise absolute errors between the reconstructed brain T2 images in Figure \ref{fig:t2qq} and the ground truth reference. The errors are scaled up by a factor of 4 for better visualization and the bright parts indicate large value. 
}
\label{fig:t2dd}
\end{figure}

\newpage
\subsubsection{Comparison with the BCD algorithm}
\begin{figure}[htbp]
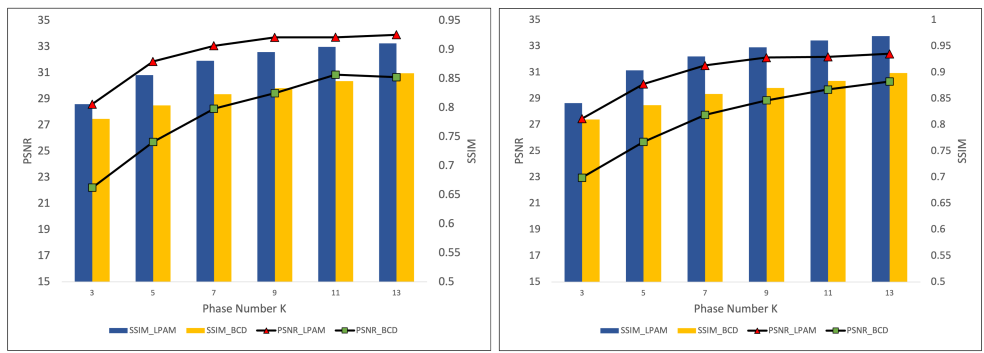

\centering
\includegraphics[width=0.4\textwidth]{T1_uv.pdf}
\includegraphics[width=0.4\textwidth]{T2uv.pdf}
\caption{Average PSNR and SSIM of the reconstruction results for images with a 10\% under-sampling ratio obtained by the LPAM algorithm and the BCD algorithm across various phase number K. Left: T1 images. Right: T2 images.}
\label{uvT1-psnr}
\end{figure}

\begin{figure}[H]
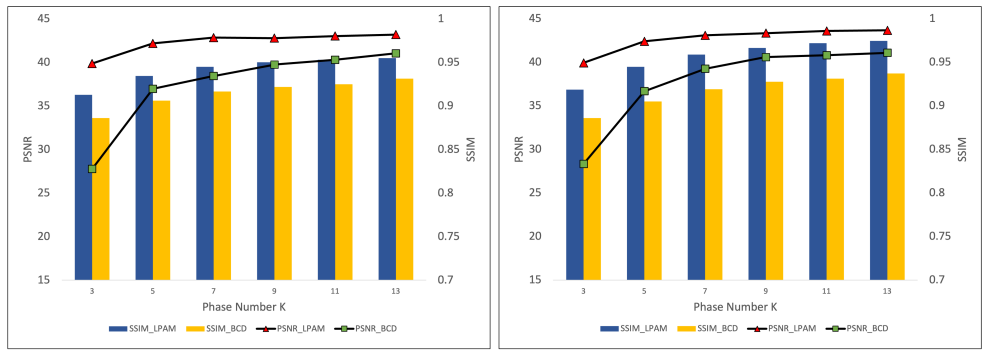

\centering
\includegraphics[width=0.4\textwidth]{T1uv20.pdf}
\includegraphics[width=0.4\textwidth]{T2uv20.pdf}
\caption{Average PSNR and SSIM of the reconstruction results for images with a 20\% under-sampling ratio obtained by the LPAM algorithm and the BCD algorithm across various phase number K. Left: T1 images. Right: T2 images.}
\label{uvT2-psnr}
\end{figure}

This experiment aims to demonstrate the effectiveness of the proposed LPAM algorithm by comparing to the standard BCD algorithm. The LPAM algorithm uses $\ubf_{k+1}$ to update the scheme to match the Res-net structure for better training, and $\vbf_{k+1}$ as a  safeguard to ensure convergence. 
{The standard BCD algorithm is equivalent to remove the iterations for updating $\ubf_{k+1}$ and only performing $\vbf_{k+1}$ in the LPAM. Although it converges for nonconvex smooth optimization, but the networks induced by BCD algorithm do not have Resnet structure.} 
We conducted experiments comparing these two methods on under-sampled T1 and T2 images using the same variational model \ref{modell}. Figure \ref{uvT1-psnr} and figure \ref{uvT2-psnr} illustrate the average PSNR and  SSIM versus phase number with respect to the radial mask with a under-sampling ratio of 10\% and 20\%. As we can tell, both the average PSNR and the average SSIM of the two algorithms ascend, indicating an increase in image quality as phase number rises. For each phase, the values of the PSNR and the SSIM of the LPAM algorithm outperforms the values of the BCD algorithm in each graph, highlighting the efficiency and accuracy of the LPAM algorithm. The compared networks were trained using Adam Optimizer implemented in Pytorch with $\beta_{1}=0.9, \beta_{2}=0.999$ and the learning rate was set at $10^{-4}.$ 

\vspace{-1cm}

\subsubsection{Comparison with the State-of-the-Art}


We compare the proposed algorithm LPAM-net with five representative methods, including the X-net\cite{do2020reconstruction}, the Joint Group Sparsity-based Network (JGSN)\cite{guo2023joint}, ReconFormer\cite{guo2023reconformer} and the Joint Cross-Attention Network (jCAN)\cite{sun2023joint}. X-net, originating from the U-net, introduces the capability to accommodate two inputs and generates two outputs and it is a non model-based method. JGSN, jCAN and ReconFormer are three model-based methods. JGSN unrolls the iterative process of the joint sparsity algorithm. jCAN deploys Vision Transformer and CNN in the image and k-space domains, respectively. ReconFormer designs Recurrent Pyamid Transformer at each architecture unit. Since the original codes for X-net and JGSN are not available, we implemented X-net and JGSN according to their original paper. 

Similar to LPAM-net, X-net and JGSN update the two under-sampled modalities jointly. ReconFormer is a single-modality reconstruction network and we reconstructed T1 and T2 images respectively for comparison. jCAN restores the under-sampled target modality with guidance from the full-sampled auxiliary modality; thus, we used full-sampled T1 images to restore under-sampled T2 images and vice versa for comparison. Since some methods incorporate transformers, we adjusted the patch size and the window size to fit our image sizes. To ensure a fair evaluation, we trained each of the state-of-the-art method for 100 epoches and the quantitative analysis of the reconstructed images is presented in the table \ref{state1} and table \ref{state2}. Compared with the other methods, the proposed LPAM-net achieves competitive performance on both T1 and T2 images with a under-sampling ratio of 20\%. 

p
\begin{table}[H]
\centering
\addtolength{\tabcolsep}{-1pt}
\caption{The comparison (mean $\pm$ standard deviation) between LPAM-net and the state-of-the-art methods for T1 images with an under-sampling ratio of 20\%}
\begin{tabular}{lccccc}
\toprule
   \textbf{Method}& 
   \textbf{PSNR} & \textbf{SSIM} & \textbf{NMSE} &
   \textbf{RMSE} & \textbf{\#Par}\\
\midrule 

X-net &  32.713$\pm$ 0.93 & 0.93$\pm$ 0.006 & 0.002$\pm$ 0.0006 & 0.024$\pm$ 0.024&42.816M\\

JGSN &  38.617$\pm$ 1.22 & 0.972$\pm$ 0.0048 & 0.0006$\pm$ 0.0002 &0.012$\pm$ 0.0016 & 21192\\

jCAN &  35.367$\pm$ 1.006 & 0.927$\pm$ 0.0087 & 0.0012$\pm$ 0.0004 & 0.017$\pm$ 0.0019 & 45.1M\\

ReconFormer& 39.124$\pm$ 1.47 & 0.977$\pm$ 0.0047 & 0.0005$\pm$ 0.0002 & 0.011$\pm$ 0.0019 & 1.1M\\
LPAM-net&  \textbf{40.66}$\pm$ 1.508 & \textbf{0.983}$\pm$ 0.004 & \textbf{0.0004}$\pm$ 0.0002 & \textbf{0.0094}$\pm$ 0.0016 & 56510\\
\bottomrule
\toprule
\end{tabular}
\label{state1}
\end{table}

\begin{table}[H]
\centering
\addtolength{\tabcolsep}{-1pt}
\caption{The comparison (mean $\pm$ standard deviation) between LPAM-net and the state-of-the-art methods for T2 images with an under-sampling ratio of 20\%}
\begin{tabular}{lccccc}
\toprule
   \textbf{Method}& 
   \textbf{PSNR} & \textbf{SSIM} & \textbf{NMSE} &
   \textbf{RMSE} & \textbf{\#Par}\\
\midrule 

X-net &  32.65$\pm$ 1.634 & 0.923$\pm$ 0.0106 & 0.0081$\pm$ 0.0033 & 0.024$\pm$ 0.005&42.816M\\

JGSN &  39.3$\pm$ 1.413 & 0.975$\pm$ 0.0046 & 0.0017$\pm$ 0.0006 &0.011$\pm$ 0.0017 & 21192\\

jCAN &  37.583$\pm$ 1.513 & 0.964$\pm$ 0.0054 & 0.0025$\pm$ 0.001 & 0.013$\pm$ 0.0023 & 45.1M\\

ReconFormer&  40.58$\pm$ 1.706 & 0.982$\pm$ 0.0051 & 0.0012$\pm$ 0.0004 & 0.0095$\pm$ 0.0018 & 1.1M\\
LPAM-net&  \textbf{42.536}$\pm$ 1.527 & \textbf{0.987}$\pm$ 0.0043 & \textbf{0.0008}$\pm$ 0.0003 & \textbf{0.0076}$\pm$ 0.0013 & 56510\\
\bottomrule
\toprule
\end{tabular}
\label{state2}
\end{table}

\section{Declarations}
The authors declare that they have no conflict of interest.

%
%

\bibliographystyle{spmpsci}      
\bibliography{reference}

\begin{thebibliography}{10}
\providecommand{\url}[1]{{#1}}
\providecommand{\urlprefix}{URL }
\expandafter\ifx\csname urlstyle\endcsname\relax
  \providecommand{\doi}[1]{DOI~\discretionary{}{}{}#1}\else
  \providecommand{\doi}{DOI~\discretionary{}{}{}\begingroup
  \urlstyle{rm}\Url}\fi

\bibitem{abolghasemi2012gradient}
Abolghasemi, V., Ferdowsi, S., Sanei, S.: A gradient-based alternating
  minimization approach for optimization of the measurement matrix in
  compressive sensing.
\newblock Signal Processing \textbf{92}(4), 999--1009 (2012)

\bibitem{aggarwal2018modl}
Aggarwal, H.K., Mani, M.P., Jacob, M.: Modl: Model-based deep learning
  architecture for inverse problems.
\newblock IEEE transactions on medical imaging \textbf{38}(2), 394--405 (2018)

\bibitem{attouch2010proximal}
Attouch, H., Bolte, J., Redont, P., Soubeyran, A.: Proximal alternating
  minimization and projection methods for nonconvex problems: An approach based
  on the kurdyka-{\l}ojasiewicz inequality.
\newblock Mathematics of operations research \textbf{35}(2), 438--457 (2010)

\bibitem{attouch2013convergence}
Attouch, H., Bolte, J., Svaiter, B.F.: Convergence of descent methods for
  semi-algebraic and tame problems: proximal algorithms, forward--backward
  splitting, and regularized gauss--seidel methods.
\newblock Mathematical Programming \textbf{137}(1-2), 91--129 (2013)

\bibitem{beck2013convergence}
Beck, A., Tetruashvili, L.: On the convergence of block coordinate descent type
  methods.
\newblock SIAM journal on Optimization \textbf{23}(4), 2037--2060 (2013)

\bibitem{bertsekas2015parallel}
Bertsekas, D., Tsitsiklis, J.: Parallel and distributed computation: numerical
  methods.
\newblock Athena Scientific (2015)

\bibitem{wanyuMiccai}
Bian, W., Zhang, Q., Ye, X., Chen, Y.: A learnable variational model for joint
  multimodal mri reconstruction and synthesis.
\newblock Lecture Notes in Computer Science \textbf{13436}, 354–364 (2022).
\newblock \doi{https://doi.org/10.1007/978-3-031-16446-0_34}

\bibitem{bolte2014proximal}
Bolte, J., Sabach, S., Teboulle, M.: Proximal alternating linearized
  minimization for nonconvex and nonsmooth problems.
\newblock Mathematical Programming \textbf{146}(1-2), 459--494 (2014)

\bibitem{charbonnier1997deterministic}
Charbonnier, P., Blanc-F{\'e}raud, L., Aubert, G., Barlaud, M.: Deterministic
  edge-preserving regularization in computed imaging.
\newblock IEEE Transactions on image processing \textbf{6}(2), 298--311 (1997)

\bibitem{LDA}
Chen, Y., Liu, H., Ye, X., Zhang, Q.: Learnable descent algorithm for nonsmooth
  nonconvex image reconstruction.
\newblock SIAM Journal on Imaging Sciences \textbf{14}(4), 1532--1564 (2021).
\newblock \doi{10.1137/20M1353368}.
\newblock \urlprefix\url{https://doi.org/10.1137/20M1353368}

\bibitem{chun2020convolutional}
Chun, I.Y., Fessler, J.A.: Convolutional analysis operator learning:
  Acceleration and convergence.
\newblock IEEE Transactions on Image Processing \textbf{29}, 2108--2122 (2020)

\bibitem{chun2020momentum}
Chun, I.Y., Huang, Z., Lim, H., Fessler, J.: Momentum-net: Fast and convergent
  iterative neural network for inverse problems.
\newblock IEEE transactions on pattern analysis and machine intelligence
  (2020)

\bibitem{chun2019bcd}
Chun, I.Y., Zheng, X., Long, Y., Fessler, J.A.: Bcd-net for low-dose ct
  reconstruction: Acceleration, convergence, and generalization.
\newblock In: Medical Image Computing and Computer Assisted
  Intervention--MICCAI 2019: 22nd International Conference, Shenzhen, China,
  October 13--17, 2019, Proceedings, Part VI 22, pp. 31--40. Springer (2019)

\bibitem{chun2018deep}
Chun, Y., Fessler, J.A.: Deep bcd-net using identical encoding-decoding cnn
  structures for iterative image recovery.
\newblock In: 2018 IEEE 13th Image, Video, and Multidimensional Signal
  Processing Workshop (IVMSP), pp. 1--5. IEEE (2018)

\bibitem{ding2023learned}
Ding, C., Zhang, Q., Wang, G., Ye, X., Chen, Y.: Learned alternating
  minimization algorithm for dual-domain sparse-view ct reconstruction.
\newblock arXiv preprint arXiv:2306.02644  (2023)

\bibitem{do2020reconstruction}
Do, W.J., Seo, S., Han, Y., Ye, J.C., Choi, S.H., Park, S.H.: Reconstruction of
  multicontrast mr images through deep learning.
\newblock Medical physics \textbf{47}(3), 983--997 (2020)

\bibitem{guo2023joint}
Guo, D., Zeng, G., Fu, H., Wang, Z., Yang, Y., Qu, X.: A joint group
  sparsity-based deep learning for multi-contrast mri reconstruction.
\newblock Journal of Magnetic Resonance \textbf{346}, 107354 (2023)

\bibitem{guo2023reconformer}
Guo, P., Mei, Y., Zhou, J., Jiang, S., Patel, V.M.: Reconformer: Accelerated
  mri reconstruction using recurrent transformer.
\newblock IEEE transactions on medical imaging  (2023)

\bibitem{hardt2014understanding}
Hardt, M.: Understanding alternating minimization for matrix completion.
\newblock In: 2014 IEEE 55th Annual Symposium on Foundations of Computer
  Science, pp. 651--660. IEEE (2014)

\bibitem{he2016deep}
He, K., Zhang, X., Ren, S., Sun, J.: Deep residual learning for image
  recognition.
\newblock In: Proceedings of the IEEE conference on computer vision and pattern
  recognition, pp. 770--778 (2016)

\bibitem{jain2013low}
Jain, P., Netrapalli, P., Sanghavi, S.: Low-rank matrix completion using
  alternating minimization.
\newblock In: Proceedings of the forty-fifth annual ACM symposium on Theory of
  computing, pp. 665--674 (2013)

\bibitem{he-1}
Kaiming~He Xiangyu~Zhang, S.R., Sun, J.: Deep residual learning for image
  recognition.
\newblock 2016 IEEE Conference on Computer Vision and Pattern Recognition p.
  770–778 (2016)

\bibitem{he-2}
Kaiming~He Xiangyu~Zhang, S.R., Sun, J.: Deep residual learning for image
  recognition.
\newblock European conference on computer vision pp. 630--645 (2016).
\newblock \doi{10.1007/978-3-319-46493-0_38}

\bibitem{kim2017deeply}
Kim, Y., Jung, H., Min, D., Sohn, K.: Deeply aggregated alternating
  minimization for image restoration.
\newblock In: Proceedings of the IEEE conference on computer vision and pattern
  recognition, pp. 6419--6427 (2017)

\bibitem{menze2014multimodal}
Menze, B.H., Jakab, A., Bauer, S., Kalpathy-Cramer, J., Farahani, K., Kirby,
  J., Burren, Y., Porz, N., Slotboom, J., Wiest, R., et~al.: The multimodal
  brain tumor image segmentation benchmark (brats).
\newblock IEEE transactions on medical imaging \textbf{34}(10), 1993--2024
  (2014)

\bibitem{nesterov2018lectures}
Nesterov, Y., et~al.: Lectures on convex optimization, vol. 137.
\newblock Springer (2018)

\bibitem{netrapalli2013phase}
Netrapalli, P., Jain, P., Sanghavi, S.: Phase retrieval using alternating
  minimization.
\newblock Advances in Neural Information Processing Systems \textbf{26} (2013)

\bibitem{palomar2010convex}
Palomar, D.P., Eldar, Y.C.: Convex optimization in signal processing and
  communications.
\newblock Cambridge university press (2010)

\bibitem{peters2009interference}
Peters, S.W., Heath, R.W.: Interference alignment via alternating minimization.
\newblock In: 2009 IEEE International Conference on Acoustics, Speech and
  Signal Processing, pp. 2445--2448. IEEE (2009)

\bibitem{pock2016inertial}
Pock, T., Sabach, S.: Inertial proximal alternating linearized minimization
  (ipalm) for nonconvex and nonsmooth problems.
\newblock SIAM Journal on Imaging Sciences \textbf{9}(4), 1756--1787 (2016).
\newblock \doi{10.1137/16M1064064}

\bibitem{polyak1964some}
Polyak, B.T.: Some methods of speeding up the convergence of iteration methods.
\newblock Ussr computational mathematics and mathematical physics
  \textbf{4}(5), 1--17 (1964)

\bibitem{sroubek2011robust}
Sroubek, F., Milanfar, P.: Robust multichannel blind deconvolution via fast
  alternating minimization.
\newblock IEEE Transactions on Image processing \textbf{21}(4), 1687--1700
  (2011)

\bibitem{sun2016deep}
Sun, J., Li, H., Xu, Z., et~al.: Deep admm-net for compressive sensing mri.
\newblock Advances in neural information processing systems \textbf{29} (2016)

\bibitem{sun2023joint}
Sun, K., Wang, Q., Shen, D.: Joint cross-attention network with deep modality
  prior for fast mri reconstruction.
\newblock IEEE Transactions on Medical Imaging  (2023)

\bibitem{ldct}
Zhang, Q., Alvandipour, M., Xia, W., Zhang, Y., Ye, X., Chen, Y.: Provably
  convergent learned inexact descent algorithm for low-dose ct reconstruction
  (2021).
\newblock \doi{10.48550/ARXIV.2104.12939}.
\newblock \urlprefix\url{https://arxiv.org/abs/2104.12939}

\end{thebibliography}

\end{document}